\documentclass[10pt]{amsart}
\usepackage{amssymb}
\usepackage{graphicx}
\usepackage{hyperref,mathrsfs}
\usepackage[all]{xy}

\newtheorem{theorem}{Theorem}[section]

\newtheorem{proposition}[theorem]{Proposition}
\newtheorem{corollary}[theorem]{Corollary}

\numberwithin{equation}{section}
\theoremstyle{definition}
\newtheorem{definition}[theorem]{Definition}
\newtheorem{example}[theorem]{Example}

\theoremstyle{remark}
\newtheorem{remark}[theorem]{Remark}

\newcommand{\abs}[1]{\lvert#1\rvert}

\newcommand{\C}{{\mathbb C}}
\newcommand{\Q}{{\mathbb Q}}
\newcommand{\R}{{\mathbb R}}
\newcommand{\Z}{{\mathbb Z}}
\newcommand{\F}{{\mathbb F}}

\newcommand{\RP}{{\mathbb R}\mathrm{P}}

\DeclareMathOperator{\SL}{{\mathrm SL}}
\DeclareMathOperator{\GL}{{\mathrm GL}}

\DeclareMathOperator{\Aut}{{\mathrm Aut}}

\DeclareMathOperator{\rank}{{\mathrm rank}}

\DeclareMathOperator{\Hom}{{\mathrm Hom}}

\DeclareMathOperator{\Tors}{{\mathrm Tors}}

\DeclareMathOperator{\Img}{{\mathrm Im}}

\begin{document}
\title{Twisted Alexander polynomial of links in the projective space}
\author{Vu Q. Huynh and Thang T. Q. Le}
\date{first version: June 12, 2005, this version: Dec 10, 2006}
\address{Faculty of Mathematics and Informatics, University of Natural Sciences, Vietnam National University, 227 Nguyen Van Cu, District 5, Ho Chi Minh City, Vietnam}
\email{hqvu@mathdep.hcmuns.edu.vn}
\urladdr{http://www.math.hcmuns.edu.vn/\~{}hqvu}
\address{School of Mathematics, Georgia Institute of Technology, Atlanta, GA 30332-0160, USA}
\email{letu@math.gatech.edu}
\urladdr{http://www.math.gatech.edu/\~{}letu}

\begin{abstract}We use Reidemeister torsion to study a twisted Alexander polynomial, as defined by Turaev, for links in the projective space. Using sign-refined torsion we derive a skein relation for a normalized form of this polynomial.
\end{abstract}

\maketitle

\section{Introduction}
  The study of polynomial invariants for links in the projective space $\RP^3$ was initiated in 1990 by Drobotukhina \cite{Drobotukhina90}. She provided a set of Reidemeister moves for links in $\RP^3$, and constructed an analogue of the Jones polynomial using Kauffman's approach involving state sum and the Kauffman bracket. Later she composed a table of links in $\RP^3$  up to six crossings, using the method of Conway's tangles \cite{Drobotukhina94}. More recently Mroczkowski \cite{Mroczkowski04} defined the Homflypt and Kauffman polynomials using an inductive argument on descending diagrams similar to the one for $S^3$.

The twisted Alexander polynomial of a link associated to a representation of the fundamental group of the link's complement  to $\GL(n;\F)$ is a generalization of the Alexander polynomial and has been studied since the early 1990s. 
In some circumstances the twisted polynomial is more powerful than the usual one: It could distinguish some pairs of knots which the usual polynomial could not, and it also provides more information on fiberedness and sliceness of knots.

For a link in $\RP^3$, the Alexander polynomial will not detect  information coming from the torsion part of the first homology group of the link's complement. We will study a version of the twisted Alexander polynomial defined by Turaev which takes the torsion part of the first homology group into account.

In his 1986 paper Turaev \cite{Turaev86}  extensively studied the Alexander polynomial using the method of Reidemeister torsion. By introducing a refinement of Reidemeister torsion -- the sign-refined torsion -- he was able to normalize the Alexander polynomial and derive a skein relation for it. Since then the sign-refined torsion has played important roles in such works as on the Casson invariant \cite{Lescop96} and the Seiberg-Witten invariant \cite{Turaev01}.

In this paper, following Turaev's method, we first identify our twisted Alexander polynomial with a corresponding  Reidemeister torsion (Theorem \ref{torsion=tAP}). Using torsion we derive a skein relation for the polynomial with a certain indeterminacy (Theorem \ref{skeintorsionthm}). Then by introducing sign-refined torsion we normalized the twisted Alexander polynomial and provide a skein relation without indeterminacies (Theorem \ref{mainthmc1}). Finally we study relationships between the twisted Alexander polynomial of a link and the  Alexander polynomial of the link's lift to $S^3$ (Theorem \ref{cover}), also using Reidemeister torsion. Although many of Turaev's arguments carry to our case,  for the sake of completeness we still provide them in details.

In our view the interest here lies primarily on the 3-dimensional nature of the method. Skein relations for link polynomial invariants are usually studied diagrammatically on two-dimensional link projections. Here we study skein relations through three-dimensional topology, using a classical yet contemporary topological invariant -- the Reidemeister torsion. 

\section{Diagrams for links in $\RP^3$ and the fundamental group}
\subsection{Diagrams} \label{diagram}
Throughout if $L$ is a link in $\RP^3$ then we let $X=\RP^3\setminus\stackrel{\circ}{N(L)}$ be its complement, where $N(L)$ is a tubular neighborhood of $L$, a collection of solid tori. We write $\pi=\pi_1(X)$ and $H=H_1(X)$.

\begin{figure}[h]
\begin{center}
\includegraphics[height=0.2\textheight]{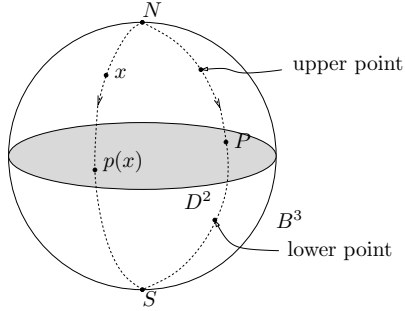}
\end{center}
\caption{Standard model.}
\label{StandardModel}
\end{figure}

We follow the terminology of Drobotukhina in \cite{Drobotukhina90}.
Consider the standard model of $\RP^3$ as a ball $B^3$ with antipodal points on the boundary sphere $\partial B^3$ identified. In this way $\RP^3=\RP^2\cup B^3$. Let $N$ and $S$ be respectively the North Pole and the South Pole of $\partial B^3$. Given a link $L$ in $\RP^3$, let $\widetilde{L}$ be its inverse image in $B^3$ under the quotient map. Isotope $L$ a bit so that $\widetilde{L}$ does not pass through $N$ or $S$. Define a projection map $p$ from $\widetilde{L}$ to the equator disk $D^2$ so that a point $x$ is mapped to the point $p(x)$ which is the intersection between the disk $D^2$ and the semicircle passing through the three points $N$, $S$ and $x$, see Fig. \ref{StandardModel}.

We can always isotope $L$ so that $\widetilde{L}$ satisfies the following conditions of general position:
\begin{enumerate}
\item $\widetilde{L}$ intersects the boundary sphere $\partial B^3$ transversally, no two points of $\widetilde{L}$ lie on the same half of a great circle joining $N$ and $S$ (i.e. $p(\widetilde{L})$ has  no double point on the boundary circle $\partial D^2$).
\item The projection $p(\widetilde{L})$ contains no cusps, no points of tangency, and no triple points.

\end{enumerate}

At each double point $P$ of $p(\widetilde{L})$, the inverse image $p^{-1}(P)$ consists of two points in $\widetilde{L}$ which are on the same semicircle joining $N$ and $S$; the one nearer to $N$ is called the upper point, the other one is called the lower point. The projection of a small arc of $\widetilde{L}$ around an upper point is called an overpass, similarly, the projection of  a small arc of $\widetilde{L}$ around a lower point is called an underpass. The projection $p(\widetilde{L})$ together with information about overpasses and underpasses is called a diagram of the link $L$. Figure \ref{2_1-knot} is an example of a diagram of a link.

\subsection{A Wirtinger-type presentation for the fundamental group}
Let $D$ be a diagram of a link $L$ (we always consider a knot as a link having one component). Choose an orientation for $D$. Label the upper arcs of $D$, each of which connecting two  underpasses, as $a_1,a_2,\dotsc,a_q,\, q\ge 0$, in arbitrary order (in case of an unknotting component which does not cross under, consider the whole component as an upper arc).  Let $2p,\ p\ge 0$ be the number of intersections  between  $D$ and the boundary circle of the projection disk. Label the intersection point counterclockwise as $b_1,b_2,\dotsc,b_{2p}$, starting from any point. To each $b_i$, associate a number $\epsilon_i$ as follows. At the point $b_i$, if $D$  is entering the boundary  then let $\epsilon_i=1$, and  let $\epsilon_i=-1$ in the other case.

Similar to the case of links in $S^3$, an application of the van Kampen gives us the following presentation for the fundamental group.
\begin{theorem}\label{presentation}
 With the above notations $\pi$ has a presentation with
generators
$a_1,a_2,\dotsc,a_q,b_1,b_2,\dotsc,b_{2p},c;$
and relations:
$$b_{p+i}=c^{-1}b_1^{\epsilon_1}b_2^{\epsilon_2}\dotsb  b_{i-1}^{\epsilon_{i-1}}b_i  b_{i-1}^{-\epsilon_{i-1}}b_{i-2}^{-\epsilon_{i-2}} \dotsb  b_1^{-\epsilon_1}  c    ,\ 1\le i\le p;$$
$$ b_1^{\epsilon_1}b_2^{\epsilon_2}\dotsb b_p^{\epsilon_p}=c^2$$
together with Wirtinger-type relations involving $a_i$'s and $b_j$'s at each crossing; and
if there is an upper arc connecting $b_i$ and $b_j$ then there is a relation $b_i=b_j$.
\end{theorem}

\begin{figure}[h]
\begin{center}
\includegraphics[height=0.2\textheight]{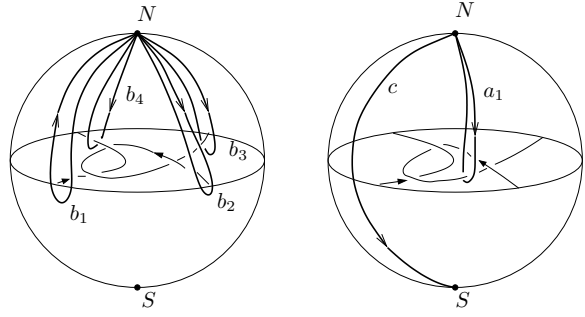}
\end{center}
\caption{Generators.}
\label{generators}
\end{figure}

\begin{figure}[h]
\begin{center}
\includegraphics[height=0.2\textheight]{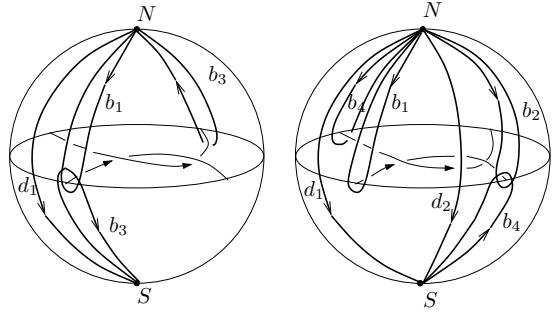}
\end{center}
\caption{Relations.}
\label{relations}
\end{figure}

\begin{remark}\label{deficiency1} In \cite{Huynh05} it is shown that if the diagram contains more
  than one crossing then a relation at a crossing can be deduced from
  the remaining relations. As a consequence if there is no affine
  unknot component then in the presentation of Theorem \ref{presentation} one may choose to omit one Wirtinger-type relation so that the number of generators is one more than the number of relations.
\end{remark}

\subsection{The first homology group}
\begin{corollary}\label{homology}
Let $L$ be a link with $v$ components, and let $S$ be the $\RP^2$ surface represented by the boundary sphere in our model. If there exists one component of $L$ whose  number of intersection points with  $S$ is odd then this component represents the non-trivial first homology class of $\RP^3$, and $H\cong\Z^v$. In the other case, $L$ represents the trivial homology class of $\RP^3$ and $H\cong\Z^v\oplus\Z_2$.
\end{corollary}
\begin{proof}
 As a result of the abelianization, the Wirtinger-type relations and the relation
$$b_{p+i}=c^{-1}b_1^{\epsilon_1}b_2^{\epsilon_2}\dotsb  b_{i-1}^{\epsilon_{i-1}}b_i  b_{i-1}^{-\epsilon_{i-1}}b_{i-2}^{-\epsilon_{i-2}} \dotsb  b_1^{-\epsilon_1}  c    , 1\le i\le p$$
 would identify all the $b_i$ and $a_j$ corresponding to the same $k$-th component of $L$ as an element $t_k\in H$, and would also identify $b_i$ and $b_{p+i}$. Thus
$$H=\langle c,t_1,t_2,\dots,t_v/ct_i=t_ic, t_it_j=t_jt_i,\prod_{i=1}^v t_i^{\delta_i}=c^2\rangle;$$
where $\delta_i$ is the sum of all $\epsilon_k,\ 0\le k\le p$, such that $b_k$ corresponds to the $i$-th component, $1\le i\le v$.

There are two cases:

{\it Case 1: All $\delta_i$ are even.} Write $\delta_i=2k_i$, $k_i\in\Z$, $1\le i\le v$. In this case $t_1^{2k_1}t_2^{2k_2}\dotsb t_v^{2k_v}=c^2$, so $(ct_1^{-k_1}t_2^{-k_2}\dotsb t_v^{-k_v})^2=1$. Let $u=ct_1^{-k_1}t_2^{-k_2}\dotsb t_v^{-k_v}$. Then $u^2=1$ and $c=ut_1^{k_1}t_2^{t_2}\dotsb t_v^{k_v}$, so
$$H=\langle t_1,t_2,\dotsc,t_v,u/t_iu=ut_i,t_it_j=t_jt_i,u^2=1\rangle\cong\Z^v\oplus\Z_2.$$

{\it Case 2: There is a $\delta_i$ that is odd.} Let $I=\{i,1\le i\le v/\delta_i=2k_i+1\}$ and $J=\{i,1\le i\le v\}\setminus I$. Let $i_0=\min\{i/i\in I\}$. Then $\prod_{i\in I}t_i^{2k_i+1}\prod_{j\in J}t_j^{2k_j}=c^2$, so $\prod_{i\in I} t_i=c^2(\prod_{i\in I}t_i^{-k_i})^2(\prod_{j\in J}t_j^{-k_j})^2=(c\prod_{1\le i\le v}t_i^{-k_i})^2$. Let $u=c\prod_{1\le i\le v}t_i^{-k_i}$. Then $\prod_{i\in I}t_i=u^2$. Since $i_0\in I$ we have $t_{i_0}\prod_{i\in I\setminus\{i_0\}}t_i=u^2$, which implies that  $t_{i_0}=u^2\prod_{i\in I\setminus\{i_0\}}t_i^{-1}$. Also $c=u\prod_{1\le i\le v}t_i^{k_i}=t^{1+2k_{i_0}}\prod_{i\in I\setminus\{i_0\}}t_i^{k_i-k_{i_0}}\prod_{j\in J}t_j^{k_j}$. So
$$H=\langle t_1,t_2,\dotsc,\hat{t_{i_0}},\dotsc,t_v,u/t_iu=ut_i,t_it_j=t_jt_i\rangle\cong\Z^v$$
(a hat over an item indicates that the item is omitted).


Consider any component $K$ of $L$. According to Poincar\'e Duality  there is a non-degenerate bilinear form $H_1(\RP^3;\Z_2)\times H_2(\RP^3;\Z_2)\rightarrow\Z_2$ where $\langle[K],[S]\rangle$ is exactly the mod $2$ intersection number between the curve $K$ and the surface $S$, 
which is $\delta_i$ mod $2$.  When $\langle[K],[S]\rangle=1$ we would have $[K]$ is non trivial in $H_1(\RP^3;\Z_2)\cong\Z_2$, and $[S]$ is non trivial in $H_2(\RP^3;\Z_2)\cong\Z_2$. On the other hand $\langle[K],[S]\rangle=0$ would imply that $[K]$ is trivial in $H_1(\RP^3;\Z_2)$.
\end{proof}

\subsubsection{Terminology}\label{classnotation} We will call a link a \emph{nontorsion link}\index{link!nontorsion} if each of its component is null-homologous (in its diagram  the number of intersection points of each component with the boundary circle is a multiple of four). The first homology group of its complement is isomorphic to $\Z^v\oplus\Z_2$. The other links are called \emph{torsion links}\index{link!torsion}. From now on we will fix the splitting of $H$ as in Corollary \ref{homology}. In this splitting if a link is nontorsion then the free part of the first homology group is generated by the meridians.

\section{Twisted Alexander polynomial}\index{twisted Alexander polynomial!Turaev's}

\subsection{The twisted homomorphism from $\Z[H]$ to $\Z[G]$}\label{twistedmap} 
Fix a splitting of $H$
as a product $H=G\times\Tors H$ of the torsion part $\Tors H=\langle u\rangle$ and
a free part $G\cong H/\Tors H$. Consider a representation (a
character) $\varphi$ from $\Tors H=\langle u\rangle$ to
$\Aut_{\C}(\C)\cong\C^*$. If $\abs{\Tors H}=1$ let  $\varphi(u)=1$; if
$\abs{\Tors H}=2$, let $\varphi(u)=-1$, i.e.
$\varphi(u)=(-1)^{\abs{\Tors H}+1}$. The map $\varphi$ then  induces a
ring homomorphism, called the \emph{twisted homomorphism} from $\Z[H]$ to $\Z[G]$ by defining
$\varphi(gu)=g\varphi(u)$. In the case $\abs{\Tors
  H}=1$, $\varphi$ is exactly the canonical projection from $\Z[H]$ to
$\Z[G]$. 
The composition of $\varphi$ and the
canonical projection $\Z[\pi]\rightarrow \Z[H]$ gives us a ring homomorphism
from $\Z[\pi]$ to $\Z[G]$. 

Let $F$ be the free group generated by the generators of $\pi$, and let $pr$ be the canonical projection $\Z[F]\rightarrow\Z[\pi]\rightarrow\Z[H]$. From now on for simplicity of notation  depending on the
context we  use the letter $\varphi$ for the twisted map above, either from
$\Z[F]$ to $\Z[G]$, or from $\Z[\pi]$ to $\Z[G]$, or from $\Z[H]$ to
$\Z[G]$.


\subsection{Twisted Alexander polynomial}
Given  a presentation $\pi=\langle x_1,\dotsc,x_n/r_1,\dotsc,r_m\rangle$ with $m=n-1$ we
construct an $m\times n$ matrix, the Alexander--Fox matrix, $[pr(\partial r_i/\partial
  x_j)]_{i,j}$, whose entries are elements of $\Z[H]$. Denote by $E(\pi)$ the ideal of $\Z[H]$
generated by the $(n-1)\times (n-1)$-minors of the Alexander--Fox
matrix. It is known that $E(\pi)$ does not depend on a presentation
of $\pi$. 

\begin{definition}\label{deftAP}The twisted Alexander polynomial of $L$ is defined as 
$\Delta^{\varphi}(L)=\gcd\varphi(E(\pi))\in\Z[G].$
\end{definition}
Note that in the unique factorization domain $\Z[G]$ the greatest common divisor is only defined up to units, which are elements of $\pm G$.

\begin{remark}\label{remarktAP}
\label{tAp=Ap}
If we replace the twisted map $\varphi$ by the canonical projection $\Z[H]\rightarrow\Z[G]$ (the torsion part of $H$ is sent to $1$) then we would get the usual Alexander polynomial $\Delta(L)$. Also, if the link $L$ is a torsion link then the twisted Alexander polynomial $\Delta^{\varphi}(L)$ is exactly the Alexander polynomial $\Delta(L)$. 
 

\end{remark}

\begin{example}
\label{2_1-knotex} Let  $K$ be the knot $2_1$ in Drobotukhina's table.
\begin{figure}[h]
\begin{center}
\includegraphics[height=0.13\textheight]{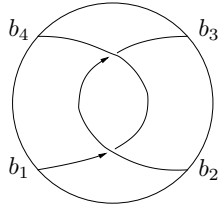}
\end{center}
\caption{The knot $2_1$.}
\label{2_1-knot}
\end{figure}
Its fundamental group has a presentation 
\begin{equation*}\begin{split}\pi&=\langle b_1,b_2,b_3,b_4,c/b_2b_1=b_4b_2=b_3b_4,b_3=c^{-1}b_1c,b_4=c^{-1}b_1^{-1}b_2b_1c,b_1^{-1}b_2^{-1}=c^2\rangle\\
&=\langle b_1,c/c=b_1c^3b_1\rangle,
\end{split}
\end{equation*}
the only relator is $r=c^{-1}b_1c^3b_1$.
Its first homology group is
$H=\langle t,c/(ct)^2=1,ct=tc\rangle$, where $t$ is the projection of the meridian $b_1$.
Let $u=ct$, then $H=\langle u,t/u^2=1,tu=ut\rangle\cong\Z\oplus\Z_2$.
The twisted homomorphism  $\varphi:\Z[\pi]\rightarrow\Z[t^{\pm 1}]$ is determined by $\varphi(b_1)=t$ and $\varphi(c)=\varphi(u)t^{-1}=-t^{-1}$.  
So $\Delta^{\varphi}_K(t)=\gcd\{\varphi(\partial r/\partial b_1),\varphi(\partial r/\partial c)\}=\gcd\{-t^{-2}(t^2-1),-t^{-1}(t-1)^2\}=t-1$. On the other hand $\Delta_K(t)=\gcd\{-t^{-2}(t^2+1),-t^{-1}(t^2+1)\}=t^2+1$.
\end{example}

\begin{example} \index{link!affine} Suppose that $L$ is an \emph{affine link}, i.e.  $L$ can be isotoped so that it is contained inside a $3$-ball in $\RP^3$, and so $L$ is a nontorsion link. Its fundamental group is  generated by $a_1,a_2,\dotsc,a_q,c$, where $q$ is the number of crossings; together with $q-1$ Wirtinger relations $r_j$ involving the $a_i$'s, and the relation $c^2=1$. 
Note that the $(q-1)\times q$ matrix $[pr(\partial r_i/\partial a_j)]_{i,j}$ is exactly the Alexander--Fox matrix of  $L$ viewed as a link in $S^3$. Then it is immediate that the twisted Alexander polynomial of $L$ is equal to $\Delta^{\varphi}(L)=\varphi(\partial c^2/\partial c)=\varphi(1+c)=0$.
On the other hand $\Delta(L)$ -- the Alexander polynomial of $L$ viewed as a link in $\RP^3$ -- will be twice the Alexander polynomial of $L$ viewed as a link in $S^3$. This supports the result that the value of the Alexander polynomial of a knot complement evaluated at $1$ is exactly the cardinality of the torsion part of the homology group (see  \cite[p. 133]{Turaev86},  \cite[p. 69]{Nicolaescu03}).
\end{example}

\section{Twisted Alexander polynomial and Reidemeister torsion}
\subsection{Background on Reidemeister torsion}
Two very readable references for this section are \cite{Milnor66} and \cite{Turaev01}.
\subsubsection{Torsion of a chain complex}
Let $\F$ be a field, $V$ be a $k$-dimensional vector space over $\F$. Suppose 
that $b=(b_1,b_2,\dotsc,b_k)$ and $c=(c_1,c_2,\dotsc,c_k)$ are two bases of $V$ then there is a non-singular $k\times k$ matrix $(a_{ij})$ such that
$c_j=\sum_{i=1}^{k}a_{ij}b_i.$
We write 
$[c/b]=\textrm{det}(a_{ij})\in\F^*$.
Two bases $b$ and $c$ are said to have the same orientation if $[b/c]>0$, and to be equivalent if $[b/c]=1$.

Let $0\rightarrow C\stackrel{\alpha}\hookrightarrow
D\stackrel{\beta}\twoheadrightarrow E\rightarrow 0$ be a short exact
sequence of vector spaces. 
Let $c=(c_1,c_2,\dotsc,c_k)$ be a basis for $C$ and
$e=(e_1,e_2,\dotsc,e_l)$ be a basis for $E$. Since $\beta$ is surjective we can
lift $e_i$ to a vector $\tilde{e}_i$ in $D$. 
Then $ce=(c_1,\dotsc,c_k,\tilde{e}_1,\dotsc,\tilde{e}_l)$ is a
basis for $D$ and its equivalence class depends not on the choice of
$\tilde{e}_i$ but only on the equivalence classes of $c$ and $e$. 

The finite chain
complex 
$(C,\partial)=(0\rightarrow
C_m\stackrel{\partial_{m}}\longrightarrow C_{m-1}\stackrel{\partial_{m-1}}
\longrightarrow
\dotsb\stackrel{\partial_2}\longrightarrow C_1\stackrel{\partial_1}\longrightarrow
C_0 \rightarrow 0)$
of finite-dimensional vector spaces over $\F$ is called \emph{acyclic} if it is exact. 
The chain is called \emph{based}
if for each $C_i$ a basis is chosen.

Assume that $(C,\partial)$ is acyclic and based with basis $c$. Choose a basis $b_i$ for
$B_i=\Img\partial_{i+1}=\ker\partial_i$.   
From the short exact sequence $0\rightarrow B_i\hookrightarrow C_i\twoheadrightarrow B_{i-1}\rightarrow 0$ we get a basis $b_ib_{i-1}$ for $C_i$. 

\begin{definition}\index{Reidemeister torsion!of chain complexes} The torsion of the acyclic and based chain complex $C$ is 
defined to be
$\tau(C)=\prod_{i=0}^{m}[b_ib_{i-1}/c_i]^{(-1)^{i+1}}\in\F^*.$
If $C$ is not acyclic then $\tau(C)$ is defined to be $0$.
\end{definition}

Note that this torsion (Turaev's version) is the inverse of 
Milnor's version. 

The torsion $\tau(C)$ depends on $c$ but does not depend on the choice of $b_i$'s.
If 
a basis $c_i'$ is used instead of $c_i$ then the torsion is multiplied with
$[c_i/c_i']^{(-1)^{i+1}}.$

\subsubsection{Torsion of a CW-complex}\label{CW-complex}
Let $X$ be a finite connected CW-complex and let $\pi=\pi_1(X)$.
The universal  cover $\widetilde{X}$ of $X$ has a canonical CW-complex
structure obtained by lifting the cells of $X$. If $\{e^k_i, 1\le i\le n_k\}$ is an ordered set of
oriented $k$-cells of $X$ and $\tilde{e}^k_i$ is any lift of $e^k_i$ 
then the ordered set $\{\tilde{e}^k_i, 1\le i\le n_k\}$ is a basis of the $\Z[\pi]$-module 
$C_i(\widetilde{X})$.

If $\Z[\pi]\stackrel{\varphi}\rightarrow\F$ is  a 
ring homomorphism then by the change of rings construction 
$\F\otimes_{\varphi}C_*(\widetilde{X})$ is a chain
complex of finite dimensional vector spaces over $\F$. If this chain complex is acyclic then we can define
its  torsion
$\tau(\F\otimes_{\varphi}C_*(\widetilde{X}))\in\F^*$.
However 
$\tau(\F\otimes_{\varphi}C_*(\widetilde{X}))$ depends on the chosen basis for $C_*(\widetilde{X})$, that is on the choices of lifting cells $\{\tilde{e}^k_i, 1\le i\le n_k\}$.
If we fix a choice of a set of lifting cells as a basis for the
$\Z[\pi]$-module $C_i(\widetilde{X})$ but change the order of the cells
in the basis then  
$\tau(\F\otimes_{\varphi}C_*(\widetilde{X}))$ is multiplied with $\pm
1$. If we change the orientations of the cells then torsion is also multiplied with $\pm 1$. If  we choose a different 
lifting cell for $e^k_i$ -- by an action $h\cdot\tilde{e}^k_i$ of a
covering  transformation  $h\in \pi$ -- then the torsion is multiplied with
$\varphi(h)^{\pm 1}$.

\begin{definition}\index{Reidemeister torsion!of CW-complexes}The Reidemeister torsion $\tau^{\varphi}(X)$ of the CW-complex
  $X$ is defined to be
the image of $\tau(\F\otimes_{\varphi}C_*(\widetilde{X}))$ under the
  quotient  map
$\F\rightarrow\F/\pm\varphi(\pi).$
\end{definition}

 Torsion is a simple homotopy invariant and a topological invariant 
of compact connected CW-complexes.
 In dimensions three or less,  where our interests are, each topological manifold has a unique piecewise-linear structure, so the torsion of a manifold can be defined.

\subsubsection{Torsion with homological bases}\index{Reidemeister torsion!with homological bases}
Here we consider the case when the chain complex is not acyclic, following \cite[p. 158]{Milnor66}.
Suppose that 
$(C,\partial)=(0\rightarrow
C_m\stackrel{\partial_{m}}\longrightarrow C_{m-1}\stackrel{\partial_{m-1}}
\longrightarrow
\dotsb\stackrel{\partial_2}\longrightarrow C_1\stackrel{\partial_1}\longrightarrow
C_0 \rightarrow 0)$  is a chain complex of based finite-dimensional vector spaces, not
necessarily acyclic.
Let $Z_i=\ker \partial_i$ and $B_i=\Img\partial_{i+1}$. Let $H_i(C)=Z_i/B_i$ and $h_i$ be its chosen basis. There is a short exact sequence
$0\rightarrow B_i\hookrightarrow Z_i\twoheadrightarrow H_i\rightarrow
0.$
This  combined with the short exact sequence
$0\rightarrow Z_i\hookrightarrow C_i\twoheadrightarrow
B_{i-1}\rightarrow 0$
show that $(b_ih_i)b_{i-1}$ is a basis for $C_i$ (and is defined up to equivalence
bases). We
can define torsion in a similar manner: $\tau(C,c,h)=\prod_{i=0}^{m}[b_ih_ib_{i-1}/c_i]^{(-1)^{i+1}}\in\F^*.$
It depends on $c$ and $h$ but does not depend on the choice of the bases $b_i$'s.

\subsubsection{Symmetry of torsion}
\label{symmetry}
Let $M$ be a compact connected orientable three-manifold.
Suppose that in the field $\F$ there is a certain ``bar'' operation so
that for all $\alpha\in\pi$,
$\overline{\varphi(\alpha)}=\varphi(\alpha^{-1})$. 
If $\partial M$ consists of tori then 
we have
$\tau^{\varphi}(M)=\overline{\tau^{\varphi}(M)}$.
For more details see \cite[p. 70]{Turaev01}.

\subsubsection{Sign-refined
    torsion}\label{signrefinedtorsion}\index{Reidemeister torsion!sign-refined}This was introduced by Turaev \cite{Turaev86} to remove the sign ambiguity of torsion.
Let $C$ be a finite based chain complex of vector spaces over $\F$. Let
$\beta_i(C)=\sum_{j=0}^{i}\dim(H_j(C))\mod 2$, $\gamma_i(C)=\sum_{j=0}^{i}\dim(C_j)\mod 2$, and
$N(C)=\sum\beta_i(C)\gamma_i(C)\mod 2$.
Let $c$ be a basis for $C$ and $h$ be a basis for $H_*(C)$. Define
$\breve{\tau}(C,c,h)=(-1)^{N(C)}\tau(C,c,h)\in\F.$
Thus $\breve{\tau}(C,c,h)$ is  $\tau(C,c,h)$ up to a sign, and they are the same when $C$ is acyclic.

A \emph{homological orientation} \index{homological orientation} for a finite CW-complex $X$ is an orientation of the finite
dimensional vector space  $\oplus_i H_i(X;\R)$. Let $h$ be a basis for
 $H_*(X;\R)$ representating a homological orientation, i.e. $h$ is a positive basis, and let $c$ be a basis for $C_*(X;\Z)$ arising from an ordered set of oriented cells of $X$, which
gives rise to a basis for  $C_*(X;\R)$. We call a lift $\tilde{c}$  of $c$ to the universal cover $\widetilde{X}$ a
\emph{fundamental family of cells}\index{fundamental family of cells}. Let 
\begin{equation}\label{signrefinedformula}\tau_0^{\varphi}(X,\tilde{c},h)=\text{sign}(\breve{\tau}(C_*(X;\R),c,h))\tau^{\varphi}(X,\tilde{c}).
\end{equation}

\begin{definition}The sign-refined torsion $\tau_0^{\varphi}(X,h)$ 
is the image of   $\tau_0^{\varphi}(X,\tilde{c},h)$ under the projection $\F\rightarrow \F/\varphi(\pi_1(X))$.
\end{definition}

This torsion has no sign ambiguity. It depends on the homological orientation but not on the order or the orientations of
the cells of $X$, since the signs of the two terms in the product change simultaneously. The choice of the number $N(C)$ is due to a change of base
formula, with it the sign-refined torsion is invariant under
simple homotopy equivalences preserving homological
orientations \cite[p. 98]{Turaev01}. 


\index{Reidemeister torsion!product formulas}
\subsubsection{Product formulas for unrefined torsion}Suppose that $0\rightarrow C'\rightarrow C\rightarrow C''\rightarrow
0$ is a short exact sequence of finite acyclic chain complexes of
vector spaces. Suppose that the bases of $C$, $C'$ and $C''$ are \emph{compatible}, in the sense that $c_i$ is equivalent to $c_i'c_i''$, then
\begin{equation}\label{productformulasforun-refinedtorsion}\tau(C)=\pm\tau(C')\tau(C'').\end{equation}

When the chains are not acyclic there is also a product formula for torsion with homological
bases. Let $h$, $h'$ and $h''$ be the bases for $H_*(C)$, $H_*(C')$ and $H_*(C'')$ respectively. The short exact sequence involving $C$, $C'$, $C''$ above gives rise to a
finite long exact sequence of homology groups
$\mathcal{H}=(\dotsb\rightarrow H_i(C')\rightarrow H_i(C)\rightarrow H_i(C'')\rightarrow
H_{i-1}(C')\rightarrow\dotsb\rightarrow H_0(C')\rightarrow H_0(C)\rightarrow H_0(C'')\rightarrow 0)$.
Since these vector spaces are based the chain $\mathcal{H}$ has a well-defined
torsion $\tau(\mathcal{H})$, which depends on $h$, $h'$ and $h''$. Suppose that the
bases of $C$, $C'$ and $C''$ are compatible, then 
\begin{equation}\tau(C,h)=\pm\tau(C',h')\tau(C'',h'')\tau(\mathcal{H}).\label{unrefinedproductformula}
\end{equation}

\subsubsection{Product formulas for sign-refined torsion}The work of keeping track of the shuffling of the bases has been done (cf. \cite[Lemma 3.4.2]{Turaev86}) in the following formula:
\begin{equation}\label{refinedproductformula}\breve{\tau}(C,c'c'',h)=(-1)^{\mu+\nu}\breve{\tau}(C',c',h')\breve{\tau}(C'',c'',h'')\tau(\mathcal{H}),\end{equation}
in other words
\begin{equation}\tau(C,c'c'',h)=(-1)^{\mu+\nu+N(C)+N(C')+N(C'')}\tau(C',c',h')\tau(C'',c'',h'')\tau(\mathcal{H}),\end{equation}
where
$\mu=\sum[(\beta_i(C)+1)(\beta_i(C')+\beta_i(C''))+\beta_{i-1}(C')\beta_{i}(C'')]\mod
2$; and
$\nu=\sum_{i=0}^m \gamma_i(C'') \gamma_{i-1}(C')\mod 2$.

\subsubsection{Homological orientations of oriented link complements}
Let $L$ be an oriented link in an oriented rational homology three-sphere $M$, and let $X$ be the link's
complement. Let $U=N(L)=\cup_{1\le i\le v} U_i$ and let $m_i$ and $l_i$ be the meridian and the
longitude of the torus boundary component $\partial U_i$.
The canonical homological orientation  of $L$ is the orientation of the vector space $H_*(X;\R)$ represented by the basis $([pt],[m_1],\dotsc,[m_v],[\partial U_1],\dotsc,[\partial U_{v-1}])$. The classes $[m_i]$ depend on the orientation of $L$ and so does the homological orientation.

\subsubsection{Reidemeister torsion associated with representations to $\SL(n;\C)$}\index{Reidemeister torsion!associated with representations}The Reidemeister torsion associated with a representation to $\text{O}(n)$ was considered by Milnor  \cite[p. 180]{Milnor66}, see also Kitano \cite{Kitano96}.
Let $X$ be a finite connected CW-complex and let $\widetilde{X}$ be its universal cover.
Let  $\rho:\pi\rightarrow
\SL(n;\C)$  be a representation of the fundamental group.
Since there is a
natural action of $\SL(n;\C)$ on $\C^n$, which is the right
multiplication of a  matrix with a vector, by using $\rho$ we can view
$\C^n$ as  a right $\Z[\pi]$-module. Thus we can form the tensor product
$C_i^{\rho}(X)=\C^n\otimes_{\Z[\pi],\rho}C_i(\widetilde{X})$, which is a vector space  over
$\C$. 
If the induced chain complex $C_*^{\rho}(X)$ is acyclic then we can define the  torsion
$\tau^{\rho}(X)=\tau(C_*^{\rho}(X))\in\C$.
Because $\rho(\pi)\subset\SL(n;\C)$ the determinant computations will
destroy some ambiguities about the choice of representing cells, so that $\tau^{\rho}(X)$ is defined up to $\pm 1$.

\begin{subsection}{Reidemeister torsion of link complements in $\RP^3$}
Let $L$ be a link in $\RP^3$. In terms of the Euler characteristic, noting  $0=\chi(\RP^3)=\chi(X\cup N(L))=\chi(X)+\chi(N(L))-\chi(X\cap N(L))$, it follows that $\chi(X)=0$. The complement $X$ is simple homotopic to a $2$-dimensional cell complex $Y$ which has one $0$-cell $\sigma^0$; $n$ $1$-cells
$\sigma^1_1,\dotsc,\sigma^1_n$; and $m$ $2$-cells 
$\sigma^2_1,\dotsc,\sigma^2_m$, where $m=n-1$. The boundary maps  are $\partial_1=0$ and
$\partial_2(\sigma^2_i)=r_i$, where $r_i$ is a word in $\sigma^1_j$, giving a presentation of the fundamental group as $\pi=\langle x_1,x_2,\dotsc,x_n/r_1,r_2,\dotsc,r_m\rangle$. This presentation is not necessarily the same as the one in Theorem \ref{presentation}, however.

Let  $\widetilde{Y}$ be the maximal abelian cover of $Y$. Consider the cellular complexes of
$\widetilde{Y}$ as modules over $\Z[H]$. 
We have a chain complex of $\Z[H]$-modules
$C_2(\widetilde{Y})\stackrel{\partial_2}\longrightarrow C_1(\widetilde{Y})\stackrel{\partial_1}\longrightarrow C_0(\widetilde{Y})\rightarrow 0$.
The boundary maps are obtained using Fox's Free Differential Calculus: 
$\partial_1(\tilde{\sigma}^1_i)=pr(x_i-1)\tilde{\sigma}^0$ and
$\partial_2(\tilde{\sigma}^2_i)=\sum_{j=1}^n pr(\frac{\partial r_i}{\partial
x_j})\tilde{\sigma}^1_j$, where the tilde sign denotes a lift of the cell to  $\widetilde{Y}$.

Denote  the quotient field $Q(\Z[G])$ of $\Z[G]$ by $\Q(G)$. 
Using the homomorphism $\varphi:\Z[H]\rightarrow \Z[G]\hookrightarrow 
\Q(G)$, construct the tensor  
$\Q(G)\otimes_{\Z[H]}C_i(\widetilde{Y})$, considered as  a vector 
space over $\Q(G)$. We have a chain complex of vector spaces over $\Q(G)$:
$$C=(\Q(G)\otimes_{\Z[H],\varphi}C_2(\widetilde{Y})\stackrel{\partial_2}\longrightarrow \Q(G)\otimes_{\Z[H],\varphi}C_1(\widetilde{Y})\stackrel{\partial_1}\longrightarrow \Q(G)\otimes_{\Z[H],\varphi}C_0(\widetilde{Y})\rightarrow 0).$$
The boundary maps are    
$[\partial_1]_i=\varphi(x_i)-1$, and  
$[\partial_2]_{i,j}=\varphi(\frac{\partial r_j}{\partial 
x_i})$, $1\le i\le n$, $1\le j\le n-1$.
Let $A=[\partial_2]^t$. 


Denote the columns of $A$ by $u_i$, $1\le i\le n$, and denote the $(n-1)\times (n-1)$ matrix obtained from $A$ by omitting the column $u_i$ by $A_i$. Since $C$ is a chain we have $0=\partial_1(\partial_2(\tilde{\sigma}^2_i))=(\sum_{j=1}^n\varphi(\frac{\partial r_i}{\partial x_j})(\varphi(x_j)-1))\tilde{\sigma}^0$, thus $\sum_{j=1}^n\varphi(\frac{\partial r_i}{\partial x_j})(\varphi(x_j)-1)=0$.
This means $\sum_{j=1}^n(\varphi(x_j)-1)u_j=0$.
For any $i>j$ we have 
\begin{equation*}
\begin{split}
(\varphi(x_j)-1)\det A_i & =\det[u_1,\dotsc,u_{j-1},(\varphi(x_j)-1)u_j,u_{j+1},\dotsc,\hat{u_i},\dotsc,u_n]\\
&=\det[u_1,\dotsc,u_{j-1},-\sum_{k\neq j}(\varphi(x_k)-1)u_k,u_{j+1},\dotsc,\hat{u_i},\dotsc,u_n]\\
&=(-1)^{i-j+1}(\varphi(x_i)-1)\det A_j.
\end{split}
\end{equation*}
Thus  for any $i$ and $j$,
\begin{equation}\label{chaincondition}
(\varphi(x_i)-1)\det A_j=\pm(\varphi(x_j)-1)\det A_i.
\end{equation}

Because $H$ has at least one free generator (Corollary \ref{homology}), the image $\varphi(\pi)$ cannot be $\{1\}$, thus there is at least one $x_i$ such that $\varphi(x_i)\neq 1$. The  property $\partial_1(\tilde{\sigma}^1_i)=(\varphi(x_i)-1)\tilde{\sigma}^0$  implies $\partial_1(\frac{1}{\varphi(x_i)-1}\tilde{\sigma}^1_i)=\tilde{\sigma}^0$, so $\partial_1$ is onto. Therefore the chain $C$ is exact if and only if $\partial_2$ is injective, which means the rank of its matrix is exactly $n-1$. Thus 
$C$ is acyclic if and only if $A$ has  a nonzero $(n-1)\times (n-1)$ minor.

The Reidemeister torsion of $C$ with respect to $\varphi$ is the torsion $\tau^{\varphi}(Y)$ of $Y$, and since torsion is a simple homotopy invariant, it is also the torsion $\tau^{\varphi}(X)$ of $X$. 

For a moment, assume that $C$ is acyclic. Take the standard bases of $\Q(G)\otimes_{\Z[H],\varphi}C_i(\widetilde{Y})$ given by $\tilde{\sigma}_j^i$ as above. A lift of $c_0=\{\tilde{\sigma_0}\}$ is $\{\frac{1}{\varphi(x_i)-1}\tilde{\sigma}^1_i\}$. Then
\begin{equation*}
\begin{split}
\tau^{\varphi}(X)&=[(\sum_{j=1}^n \varphi(\frac{\partial r_1}{\partial x_j})\tilde{\sigma}^1_j,\dotsc,\sum_{j=1}^n \varphi(\frac{\partial r_{n-1}}{\partial x_j})\tilde{\sigma}^1_j, \frac{1}{\varphi(x_i)-1}\tilde{\sigma}^1_i)/(\tilde{\sigma}_1^1,\dotsc,\tilde{\sigma}_n^1)]\\
&=\frac{(-1)^{i+n}}{\varphi(x_i)-1}\det A_i.
\end{split}
\end{equation*}

Thus if $\varphi(x_i)\neq 1$ then
$\tau^{\varphi}(X)=\pm \det A_i/(\varphi(x_i)-1).$
By Eq. (\ref{chaincondition}) if $\varphi(x_j)=1$ then $\det(A_j)=0$, 
hence the following formula is correct for all $i$, whether $C$ is acyclic or not:
\begin{equation}\label{torsion2}
(\varphi(x_i)-1)\tau^{\varphi}(X)=\pm \det A_i\in\Q(G)/\pm G.
\end{equation}

\begin{theorem}\label{torsion=tAP}\index{link!nontorsion}
The Reidemeister torsion and the twisted Alexander polynomial of the  complement of a nontorsion link  are the same.
\end{theorem}
\begin{proof}
According to Definition \ref{deftAP} and Formula (\ref{torsion2}), we have 
$$\Delta^{\varphi}(X)=\gcd\{\det A_1,\dotsc,\det A_n\} = \gcd\{(\varphi(x_1)-1)\tau^{\varphi}(X),\dotsc,(\varphi(x_n)-1)\tau^{\varphi}(X)\}.$$
Thus we will get $\Delta^{\varphi}(X)=\tau^{\varphi}(X)$ immediately from the following claim.

{\it Claim.} $\gcd\{\varphi(x_1)-1,\varphi(x_2)-1,\dotsc,\varphi(x_n)-1\}=1\in\Z[G]/\pm G$.

To prove the claim we consider two cases.

{\it Case 1: $L$ has one component.} In this case $H=\langle t,u/tu=ut,u^2=1\rangle$, $pr(x_i)=t^{m_i}u^{n_i}$, and   $\varphi(x_i)=t^{m_i}(-1)^{n_i}$. Let $d=\gcd\{\varphi(x_1)-1,\varphi(x_2)-1,\dotsc,\varphi(x_n)-1\}\in\Z[t^{\pm 1}]$. The following two identities :
$$(t^{m_i}(-1)^{n_i}-1)+t^{m_i}(-1)^{n_i}(t^{m_j}(-1)^{n_j}-1)=t^{m_i+m_j}(-1)^{n_i+n_j}-1,$$
$$(t^{m_i}(-1)^{n_i}-1)-t^{m_i-m_j}(-1)^{n_i-n_j}(t^{m_j}(-1)^{n_j}-1)=t^{m_i-m_j}(-1)^{n_i-n_j}-1,$$
(compare \cite[p.  117]{Lickorish97}) imply that $d|(t^{\sum_{i=1}^n \alpha_i m_i}(-1)^{\sum_{i=1}^n \alpha_i n_i}-1)$ for any $\alpha_i\in\Z$.

Since $t\in pr(\pi)$ there are $\alpha_i\in \Z$ such that $t=\prod_{i=1}^n pr(x_i^{\alpha_i})=t^{\sum_{i=1}^n \alpha_i m_i}u^{\sum_{i=1}^n \alpha_i}$, which implies $\sum_{i=1}^n \alpha_i m_i=1$ and $\sum_{i=1}^n \alpha_i$ is even. Thus $d|(t-1)$, hence either $d=1$ or $d=t-1$, up to $\pm t^k,k\in\Z$. Since $u\in pr(\pi)$ there is at least an $i_0$ such that $n_{i_0}$ is odd, so that $\varphi(x_{i_0})-1=-t^{m_{i_0}}-1$. 
Since $\gcd\{t-1,-t^{m_{i_0}}-1\}=1$, we conclude that $d=1$.

{\it Case 2: $L$ has at least two components.} Let $v\ge 2$ be the number of components. Now $pr(x_i) = t_1^{m_i^1} t_2^{m_i^2} \dotsb t_v^{m_i^v} u^{n_i}$ and $\varphi(x_i)=t_1^{m_i^1}t_2^{m_i^2}\dotsb t_v^{m_i^v}(-1)^{n_i}$. Letting $t_2=t_3=\dotsb=t_v=1$ and applying the argument in 
Case 1 to $t_1$ we have the result.
\end{proof}



\begin{theorem}\label{torsionlinkAp=R}
If $L$ is a torsion knot and $t$ is the generator of the first homology group then
$\tau_L^{\varphi}(t)=\Delta_L^{\varphi}(t)/(t-1)\in\Z[t^{\pm 1},(t-1)^{-1}]$. If $L$  is  a torsion link with a least two components then 
the Reidemeister torsion and the twisted Alexander polynomial are the same. 
\end{theorem}

\begin{proof} The proof is similar to the proof of Theorem \ref{torsion=tAP}. 

{\it Case 1: $L$ has one component.} In this case $H=\langle t\rangle$ and  $\varphi(x_i)=t^{m_i}$. Using the two identities:
$t^{m_i}+t^{m_i}(t^{m_j}-1)=t^{m_i+m_j}-1$,
and $(t^{m_i}-1)-t^{m_i-m_j}(t^{m_j}-1)=t^{m_i-m_j}-1$,
we get $\gcd\{\varphi(x_i)-1,1\le i\le n\}=t-1$,
 thus $\Delta^{\varphi}(X)=(t-1)\tau^{\varphi}(X)$.

{\it Case 2: $L$ has at least two components.} Now $H$ is generated by $t_1,t_2,\dotsc,t_v$; $v\ge 2$, and $\varphi(x_i)=t_1^{m_i^1}t_2^{m_i^2}\dotsb t_v^{m_i^v}$. By subsequently letting $t_j=1$ for all $j\neq i$ and applying the argument in Case 1 to $t_i$ we obtain $\gcd\{\varphi(x_1)-1,\varphi(x_2)-1,\dotsc,\varphi(x_n)-1\}=\gcd\{t_1-1,t_2-1,\dotsc,t_v-1\}=1$, hence $\Delta^{\varphi}(X)=\tau^{\varphi}(X)$.
\end{proof}

\begin{remark}\label{Ap=R}
With a virtually identical proof, the statement of Theorem \ref{torsionlinkAp=R} is true for all links if we replace the twisted map $\varphi$ by the canonical projection $\Z[H]\rightarrow\Z[G]$ and replace the twisted Alexander polynomial by the Alexander polynomial.
\end{remark}
\end{subsection}

\begin{subsection}{Comparison with other twisted Alexander polynomials} Among the first people who studied twisted Alexander polynomials were Lin \cite{Lin01}, 
Wada \cite{Wada94}, Kitano \cite{Kitano96}, Kirk-Livingston \cite{Kirk-Livingston99a}.
Except for  Lin's construction which used Seifert surfaces, other constructions were  based on that of Wada.  The twisted Alexander polynomial in the form considered in this paper was defined first by Turaev in \cite{Turaev02a} and was discussed further in \cite[p. 27]{Turaev02}. It receives attention recently in \cite{Heusener-Porti05}. We outline Wada's construction to show its relationship with our polynomial\index{twisted Alexander polynomial!Wada's}. 

Suppose $\pi=\langle x_1,\dotsc,x_m/ r_1,\dotsc,r_{m-1}\rangle$ is a presentation of deficiency one, and let $\rho:\pi\rightarrow\GL(n;\F)$ be a representation of $\pi$. Let $\alpha:\pi\rightarrow G \cong \langle t_1,t_2,\dotsc,t_v/t_i t_j=t_j t_i\rangle\cong\Z^v$ be a surjective group homomorphism.  Define a ring homomorphism $\phi:\Z[\pi]\rightarrow M(n;\F[G])$ by letting $\phi(x)=\alpha(x)\rho(x)$ for $x\in\pi$ then extend linearly
(or equivalently one first extends $\alpha$ linearly to  a ring homomorphism $\tilde{\alpha}:\Z[\pi]\rightarrow\Z[G]$, and extends $\rho$ linearly to a ring homomorphism $\tilde{\rho}:\Z[\pi]\rightarrow M(n;\F)$, then let $\phi=\tilde{\alpha}\otimes\tilde{\rho}$).

Consider the $(m-1)\times m$ matrix $M$ whose the $(i,j)$ entry is $\phi(\partial r_i/\partial x_j)\in M(n;\F[G])$. Let $M_j$ be the $(m-1)\times (m-1)$  matrix obtained from $M$ by removing the $j$-th column. View $M_j$ as an $n(m-1)\times n(m-1)$ matrix whose entries are in $\F[G]$. Supposing that $\phi(x_j)\neq I$, we define the twisted Alexander polynomial as
$\Delta^{\rho}(X)= \det M_j/\det \phi(1-x_j)\in\F(G)=Q(\F[G]).$
Wada proved that this polynomial is independent of the choice of $j$ and the choice of a presentation of $\pi$, and is defined up to a factor in $\pm G$.

Let us compare Wada's polynomial with that of Turaev. Fix a splitting $H=G\times \Tors H$. Suppose that $\varphi\in \Hom(\Tors H,\C^*)$ is given. Let $\alpha$ be as above, and $\rho$ be the composition of the maps $\pi\rightarrow H\stackrel{\beta}\rightarrow\Tors H\stackrel{\varphi}\rightarrow\{\pm 1\}\subset\GL(1;\C)$; here the first arrow is the canonical projection map, and $\beta$ maps an element $gh\in H$ where $g\in G$ and $h\in \Tors H$ to $h$. Then $\phi=\tilde{\alpha}\otimes\tilde{\rho}$ is exactly the twisted map in Section \ref{twistedmap}. Thus , in view of Formula (\ref{torsion2}) $\Delta^{\rho}(X)$ here is exactly the torsion $\tau^{\varphi}(X)$, and its relationships with Turaev's polynomial are provided in Theorems \ref{torsion=tAP} and \ref{torsionlinkAp=R}. Unlike the general case Turaev's invariant is still abelian.

\begin{remark}
Milnor   proved in \cite{Milnor62} the identification between Alexander polynomial and Reidemeister torsion for knot complements in $S^3$.  Kitano \cite{Kitano96} proved the identification between Wada's twisted Alexander polynomial and Reidemeister torsion, also for knot complements in $S^3$. Kirk--Livingston \cite{Kirk-Livingston99a} generalized this result to CW-complex, but considered only a one variable twisted Alexander polynomial associated with an infinite cyclic cover of the complex. Turaev \cite[p. 28]{Turaev02} has also studied this problem. The elementary proof above of Theorem \ref{torsion=tAP} is  close to Milnor's original one.
\end{remark}
\end{subsection}

\begin{section}{A skein relation for the twisted Alexander polynomial}\label{firstskein}
\begin{subsection}{The one variable twisted Alexander polynomial} 
Let $L$ be a link with $v$ components, and let $\varphi'$ be the composition of $\varphi$ with the canonical
  projection from $\Z[t_1^{\pm 1},\dotsc,t_v^{\pm 1}]$ to $\Z[t^{\pm
      1}]$. The twisted Alexander polynomial $\Delta^{\varphi}(X)$ of
  the complement $X$ is
  a polynomial in $v$ variables $t_1,t_2,\dotsc,t_v$. 
  The one variable
  polynomial $\Delta^{\varphi'}(X)$ is obtained from
  $\Delta^{\varphi}(X)$ by replacing $\varphi$ by $\varphi'$. 

We write $\Q(t)=Q(\Z[t,t^{-1}])$, 
 $\Delta^{\varphi'}(X)$ as
  $\Delta_L^{\varphi'}(t)$, and $\tau^{\varphi'}(X)$ as $\tau_L^{\varphi'}(t)$. 

The proofs of the following two theorems are identical to the proofs for the cases of knots of Theorems \ref{torsion=tAP} and \ref{torsionlinkAp=R}. 

\begin{theorem}\label{otAPnontorsion}\index{link!nontorsion}
If $L$  is a nontorsion link then the Reidemeister torsion $\tau_L^{\varphi'}(t)$ and the one variable twisted Alexander polynomial $\Delta_L^{\varphi'}(t)$ are the same.
\end{theorem}

\begin{theorem}\label{otAPtorsion}
If $L$  is  a torsion link then the Reidemeister torsion and the
one variable twisted Alexander polynomial are
related by the formula
$\tau_L^{\varphi'}(t)=\Delta_L^{\varphi'}(t)/(t-1)\in\Z[t^{\pm 1},(t-1)^{-1}]$.
\end{theorem}


As a consequence of the symmetry of  torsion (Section \ref{symmetry}), we have:
\begin{theorem}\label{symmetrytorsion}
The Reidemeister torsion $\tau_L^{\varphi'}(t)$ is symmetric, that is $\tau^{\varphi'}_L(t^{-1})=\tau^{\varphi'}_L(t)$ up to $\pm t^n,n\in\Z$, as elements in $\Q(t)$.
\end{theorem}

From this we derive:
\begin{theorem}
The one variable twisted Alexander polynomial is symmetric, that is $\Delta_L^{\varphi'}(t^{-1})=\Delta_L^{\varphi'}(t)$ up to $\pm t^n,n\in\Z$, as elements in $\Z[t^{\pm 1}]$.
\end{theorem}

\begin{proof}
If $L$ is a nontorsion link then according to Theorem \ref{otAPnontorsion}, $\Delta_L^{\varphi'}(t^{-1})=\tau_L^{\varphi'}(t^{-1})=\tau_L^{\varphi'}(t)=\Delta_L^{\varphi'}(t^{-1})$. If $L$ is a torsion link then according to Theorem \ref{otAPtorsion}, $\Delta_L^{\varphi'}(t^{-1})=(t^{-1}-1)\tau_L^{\varphi'}(t^{-1})=t^{-1}(1-t)\tau_L^{\varphi'}(t)=\Delta_L^{\varphi'}(t)$, up to $\pm t^n,n\in\Z$.
\end{proof}

\end{subsection}

\begin{subsection}{A skein relation for torsion with indeterminacies}\label{indeterminacies}
Let $L$ be an oriented link. Consider a crossing of $L$. Let $B$ be an open $3$-ball that encloses
this crossing and intersects $L$ at four points. Let
$V=\RP^3\setminus (B\cup \stackrel{\circ}{N(L)})$, see Fig. \ref{skein2}.

\begin{figure}[h]
\begin{center}
\includegraphics[width=0.8\textwidth]{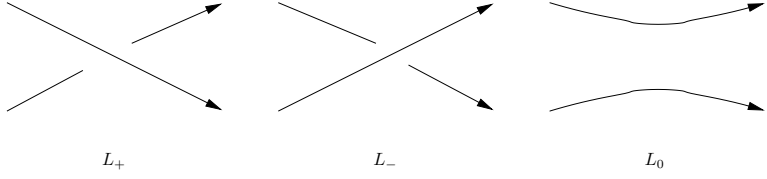}
\end{center}
\caption{The links $L_+$, $L_-$, $L_0$ are identical except at one crossing.}
\label{smoothing}
\end{figure}

\begin{figure}[h]
\begin{center}
\includegraphics[height=0.17\textheight]{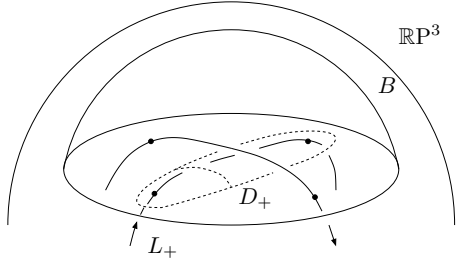}
\end{center}
\caption{The disk $D_+$ at a crossing of the link $L_+$.}
\label{skein2}
\end{figure}

 Take a triangulation of
 $V$. There is a deformation retraction (actually a simple homotopy) of the complement of
 $L_{\alpha}$, $\alpha\in\{+,-,0\}$, onto $X_{\alpha}=V\cup D_{\alpha}$, where
 $D_{\alpha}$ is a disk glued to $\partial V$ along a simple loop
 $\partial D_{\alpha}$ circling two intersection points of $B$ and $L_{\alpha}$
 as in Fig. \ref{skein3} such that $V\cup D_{\alpha}$ has a cell
 decomposition consists of the cells of $V$ plus the disk
 $D_{\alpha}$. We can assume that the loops $\partial D_{\alpha}$'s have
 a common point.
\begin{figure}[h]
\begin{center}
\includegraphics[height=0.17\textheight]{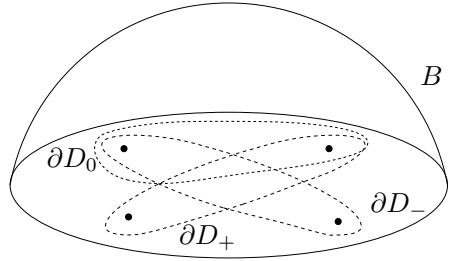}
\end{center}
\caption{The curves $\partial D_{\alpha}$'s.}
\label{skein3}
\end{figure}

\subsubsection{Smoothing of crossings and torsion classes}
When the smoothing operation is done at a particular crossing the link $L_0$ may no longer be in the same torsion class with $L_+$ and $L_-$ (recall Section \ref{classnotation}). The three links $L_+$, $L_-$ and $L_0$ are in the same torsion class in the following cases: 
\begin{enumerate}
\item There is one component of $L_+$ which is not involved at the crossing that is not null-homologous. In this case
$L_+$, $L_-$ and $L_0$ are all torsion links.

\item The two strands of $L_+$ at the crossing come from one component, and after smoothing all components are null-homologous. In this case $L_+$, $L_-$ and $L_0$ are all nontorsion links (cf. Fig. \ref{1component}).

\item The two strands of $L_+$ at the crossing come from two different components, and before smoothing all components are null-homologous. In this case $L_+$, $L_-$ and $L_0$ are all nontorsion links (cf. Fig. \ref{2components}).
\end{enumerate}

In what follows we will need the condition that the links $L_{\alpha}$'s belong to the same torsion class, so that $\Tors H_1(X_{\alpha})$'s are the same. Therefore throughout the rest of Section \ref{firstskein} we will assume that this condition is satisfied at the crossing under consideration.

\begin{subsubsection}{The chain complexes $C_{\alpha}$ and $C$}\label{C-Calpha}
Fix an $\alpha\in\{+,-,0\}$. Let $\widetilde{X}_{\alpha}$ be the $D=\Z\times\Tors H_1(X_{\alpha})$ cover of $X_{\alpha}$ corresponding to the kernel of the map $proj_{\alpha}:\pi_1(X_{\alpha})\rightarrow H_1(X_{\alpha})\rightarrow G\times\Tors H_1(X_{\alpha})\rightarrow \{t^m:m\in\Z\}\times\Tors H_1(X_{\alpha})$. Let $\widetilde{V}$ be the inverse image of $V$ under the covering map. The triangulation of $V$ induces a CW-complex structure on $\widetilde{V}$.

Under the condition that $L_{\alpha}$'s are in the same torsion class we can construct $\widetilde{X}_{\alpha}$ in a different way as follows. Take $\widetilde{V}$ to be the cover of $V$ corresponding to the kernel of the projection $\pi_1(V)\rightarrow H_1(V)\rightarrow G\times\Tors H_1(V)\rightarrow \{t^m:m\in\Z\}\times\Tors H_1(V)$. Noting that $\Tors H_1(V)=\Tors H_1(X_{\alpha})$ for $\alpha=+,-,0$, we construct $\widetilde{X}_{\alpha}$ from $\widetilde{V}$ by gluing $\lvert\Z\times\Tors H_1(X_{\alpha})\rvert$ copies of $D_{\alpha}$ along the lifts of $\partial D_{\alpha}\subset V$.

Consider the ring homomorphism $\varphi':
\Z[\{t^m:m\in\Z\}\times\Tors H_1(X_{\alpha})]\rightarrow\Z[t^{\pm1}]\hookrightarrow\Q(t)$, which does not depend on $\alpha$.
Let $C_{\alpha}=\Q(t)\otimes_{\Z[\Z\times\Tors H_1(X_{\alpha})],\varphi'}
C_*(\widetilde{X}_{\alpha};\Z)$ and let $C=\Q(t)\otimes_{\Z[\Z\times\Tors H_1(X_{\alpha})],\varphi'}
C_*(\widetilde{V};\Z)$, both considered as  chain complexes of
$\Q(t)$-vector spaces. Note that $C$ does not depend on $\alpha$.
\end{subsubsection}

\begin{subsubsection}{Relations among $\tau(C_{\alpha})$'s}
Since $C_i(\widetilde{V};\Z)\hookrightarrow
C_i(\widetilde{X}_{\alpha};\Z)$ is an inclusion, the induced map $\Q(t)\otimes_{\varphi'}
C_*(\widetilde{V};\Z)\hookrightarrow \Q(t)\otimes_{\varphi'}
C_*(\widetilde{X}_{\alpha};\Z)$ is injective, and we have the short exact sequence of chain
complexes of $\Q(t)$-vector spaces
\begin{equation}\label{pair}
0\rightarrow C\rightarrow C_{\alpha}\rightarrow C_{\alpha}/C\rightarrow 0.
\end{equation}

Choose a fundamental family of  cells 
for $\widetilde{V}$ providing a
basis for the chain $C$. A fundamental family of cells of
$\widetilde{X}_{\alpha}$ is obtained from the one of $\widetilde{V}$ by
adding a lift of $D_{\alpha}$. We can choose these lifts $D_{\alpha}$ so that the
loops $\partial\widetilde{D}_{\alpha}$ have a common point in $\widetilde{V}$, which is a
lift of the common  point of $\partial D_{\alpha}$ in $V$. 
Recalling that $X_{\alpha}$ is the result of gluing the disk $D_{\alpha}$ to $V$, 
we observe that only the second homology group of $C_{\alpha}/C$ is non-trivial, and 
the torsion of $C_{\alpha}/C$ with homology bases is $\tau(C_{\alpha}/C,h)=1$ up to a sign.

Suppose that the chain complex $C_{\alpha}$ is acyclic. The product formula for torsion (\ref{unrefinedproductformula}) applied to the short exact sequence (\ref{pair}) gives:
\begin{equation}
\tau(C_{\alpha})=\pm\tau(C,h)\tau(C_{\alpha}/C,h)\tau(\mathcal{H}_{\alpha})=\pm\tau(C,h)\tau(\mathcal{H}_{\alpha}),\label{torsionofC_i}
\end{equation}
where $\mathcal{H}_{\alpha}$ denotes the long exact homological sequence of the pair $(C_{\alpha},C)$, with a chosen basis:
$\mathcal{H}_{\alpha}=(\dotsb\rightarrow H_i(C)\rightarrow H_i(C_{\alpha})\rightarrow H_i(C_{\alpha}/C)\rightarrow H_{i-1}(C)\rightarrow\dotsb\rightarrow H_0(C)\rightarrow H_0(C_{\alpha})\rightarrow H_0(C_{\alpha}/C)\rightarrow 0)$.
Since $C_{\alpha}$ is exact the sequence $\mathcal{H}_{\alpha}$ is reduced to 
$0\rightarrow H_2(C_{\alpha}/C)\stackrel{\partial}\rightarrow H_1(C)\rightarrow 0\rightarrow 0\rightarrow 0\rightarrow 0\rightarrow 0\rightarrow 0$,
so $H_1(C)\cong H_2(C_{\alpha}/C)\cong\Q(t)$ and $\tau(\mathcal{H}_{\alpha})=\det(\partial)$.
Let $y$ be the chosen basis of the one-dimensional $\Q(t)$-vector space $H_1(C)$. Then $\partial[\widetilde{D}_{\alpha}]=[\partial\widetilde{D}_{\alpha}]=\gamma_{\alpha} y$ for some $\gamma_{\alpha}\in\Q(t)$.
Formula (\ref{torsionofC_i}) now gives
\begin{equation}\label{skeineqn2}
\tau(C_{\alpha})=\pm\gamma_{\alpha}\tau(C,h).
\end{equation}
As can be seen from the proof the product formula in \cite[p. 160]{Milnor66} or from the corresponding formula for sign-refined torsion (\ref{refinedproductformula}), the sign $\pm$ in (\ref{skeineqn2}) above depends only on the ranks of the vector spaces in the chains $C_{\alpha}$, $C$ and $\mathcal{H}_{\alpha}$, thus does not depend on $\alpha$ (to the extent that $C_{\alpha}$ is assumed to be acyclic).

Under the assumption that there is at least one $\alpha_0\in\{+,-,0\}$ such that $C_{\alpha_0}$ is acyclic, we show that  (\ref{skeineqn2}) above still holds  when $C_{\alpha}$ is not acyclic. When $C_{\alpha}$ is not acyclic, by definition $\tau(C_{\alpha})=0$. We will show that $\gamma_{\alpha}$ is zero, i.e. the boundary map  $\partial: H_2(C_{\alpha}/C)\rightarrow H_1(C)$ is  zero. Suppose the contrary, $\gamma_{\alpha}\neq 0$. Because $H_1(C)\cong H_2(C_{\alpha_0}/C)\cong\Q(t)\cong H_2(C_{\alpha}/C)$, if $\partial$ is not zero it must be a bijection. The long exact sequence $\mathcal{H}_{\alpha}$ shows that  $H_1(C_{\alpha})=0$. Note that $\rank (\Q(t)\otimes_{\Z[\Z\times\Tors H_1(X_{\alpha})],\varphi'} C_i(\widetilde{X}_{\alpha},\Z))$ is exactly the number of $i$-cells of $X_{\alpha}$. This implies that $0=\chi(X_{\alpha})=\chi(C_{\alpha})=\rank(H_0(C_{\alpha}))+\rank(H_2(C_{\alpha}))$. Thus $H_0(C_{\alpha})=H_1(C_{\alpha})=H_2(C_{\alpha})=0$ i.e. $C_{\alpha}$ is acyclic, a contradiction.
\end{subsubsection}

\begin{subsubsection}{Relations among $\gamma_{\alpha}$'s}In view of  (\ref{skeineqn2}) to further study  relations among $\tau(C_{\alpha})$'s we now try to find a relation among $\gamma_{\alpha}$'s. Recall that $\gamma_{\alpha}\in H_1(C)$ is represented  by the loop $\partial\widetilde{D}_{\alpha}$.
Let $a, b, c,d$ be simple meridian loops with a common base point, circling the four intersection points between $L$ and $B$ as in Fig. \ref{skein4}.
\begin{figure}[h]
\begin{center}
\includegraphics[height=0.17\textheight]{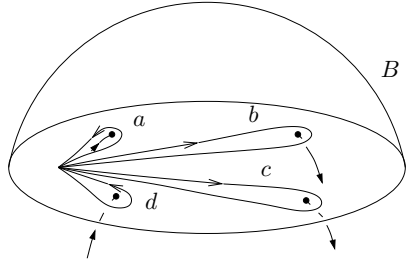}
\end{center}
\caption{The loops $a$, $b$, $c$ and $d$.}
\label{skein4}
\end{figure}
The boundary of the disks $D_{\alpha}$'s are: $\partial
D_+=bd^{-1}$, $\partial D_-=a^{-1}c$, and $\partial D_0=a^{-1}b$. 
Under the  the map $proj_{\alpha}$ in Section \ref{C-Calpha}, all of $a,b,c,d$ are projected to $t$. Noting that $a^{-1}bcd^{-1}=1$, we have $\tilde{d}=-t^{-1}\tilde{a}+t^{-1}\tilde{b}+\tilde{c}$. Hence $\gamma_0y=\widetilde{a^{-1}b}=t^{-1}(-\tilde{a}+\tilde{b})$, $\gamma_-y=\widetilde{a^{-1}c}=t^{-1}(-\tilde{a}+\tilde{c})$, and $\gamma_+y=\widetilde{bd^{-1}}=\tilde{b}-\tilde{c}+t^{-1}(\tilde{a}-\tilde{b})=(t-1)\gamma_0y-t\gamma_-y$.
So in $\Q(t)$:
\begin{equation}\label{gammas}\gamma_+ +(1-t)\gamma_0+t\gamma_-=0.\end{equation}
\end{subsubsection}

Formulas (\ref{skeineqn2}) and (\ref{gammas}) now give us, under the assumption that there is at least one $\alpha_0\in\{+,-,0\}$ such that $C_{\alpha_0}$ is acyclic, the formula $\tau(C_+)+(1-t)\tau(C_0)+t\tau(C_-)=0$. But this formula is also trivially correct when none of the $C_{\alpha}$ are acyclic, since in that case all three torsions are zero. Thus we obtain the following theorem:
\begin{theorem}\label{skeintorsionthm}If $L_+$, $L_-$ and $L_0$ belong to the same torsion class then
\begin{equation} \label{skeintorsion}\tau_{L_+}^{\varphi'}(t)+(1-t)\tau_{L_0}^{\varphi'}(t)+t\tau_{L_-}^{\varphi'}(t)=0.
\end{equation}
\end{theorem}

\end{subsection}
\begin{subsection}{Sign-refined torsion and a normalized one variable twisted Alexander function}
\begin{subsubsection}{A skein relation for sign-refined torsion}We consider sign-refined torsion, see Section \ref{signrefinedtorsion}.
In all that follow the bases  for the chain complexes are induced from the  triangulations of the spaces as previously mentioned at the beginning of Section \ref{indeterminacies}. There are two cases:

{\it Case 1: The two strands of $L_+$ at the crossing come from the same
  component.} See Fig. \ref{1component}. 
\begin{figure}[h]
\begin{center}
\includegraphics[height=0.17\textheight]{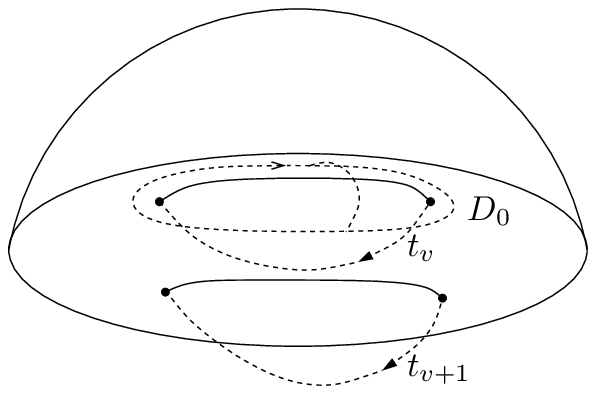}
\end{center}
\caption{Case 1.}
\label{1component}
\end{figure}
Suppose that the crossing involves the $v$-th component of
$L_+$. The bases $h_{\alpha}$ for $H_*(X_{\alpha};\R)$, $\alpha=+,-$  consist of
$[pt]$, $t_1,\dotsc,t_v$, $q_1,\dotsc,q_{v-1}$, where $q_i$
represents the $i$-th boundary component of $L_+$ and $t_i$ represent the (oriented) meridian of this component. The basis
for $H_*(X_{0};\R)$ consists of $[pt]$, $t_1,\dotsc,t_{v+1}$, $q_1,\dotsc,q_{v}$. The basis $h_0$ for $H_*(V;\R)$ consists of $[pt]$, $t_1,\dotsc,t_{v+1}$, $q_1,\dotsc,q_{v-1}$.

We want to compare the terms  $\breve{\tau}(C_*(X_{\alpha},\R),c_{\alpha},h_{\alpha}))$. Consider the short exact sequence of chain complexes:
$0\rightarrow C_*(V;\R)\rightarrow C_*(X_{\alpha};\R)\rightarrow
C_*(X_{\alpha},V;\R)\rightarrow 0.$
Applying the product formula for sign-refined torsion (\ref{refinedproductformula}) we obtain 
$$\breve{\tau}(C_*(X_{\alpha};\R))=(-1)^{\mu_{\alpha}+\nu}\breve{\tau}(C_*(V;\R))\breve{\tau}(C_*(X_{\alpha},V;\R))\tau(\mathcal{H}_{\alpha}),$$
 where $\mathcal{H}_{\alpha}$ is the long exact homological sequence of the pair $(X_{\alpha},V)$ with real coefficients, and
\begin{multline*}\mu_{\alpha}=\sum[(\beta_i(C_*(X_{\alpha};\R))+1)(\beta_i(C_*(V;\R))+\beta_i(C_*(X_{\alpha},V;\R)))+\\
+\beta_{i-1}(C_*(V;\R))\beta_{i}(C_*(X_{\alpha},V;\R))]\mod
2
\end{multline*}
and
$\nu= \sum_{i=0}^m\gamma_i(C_*(X_{\alpha},V;\R))\gamma_{i-1}(C_*(V;\R))\mod 2$.
Notice that $\nu$ does not depend on $\alpha$.

Since the term $\breve{\tau}(C_*(V;\R))\breve{\tau}(C_*(X_{\alpha},V;\R))$ does not depend on $\alpha$ we only need to compare the terms $(-1)^{\mu_{\alpha}}\text{sign}(\tau(\mathcal{H}_{\alpha}))$.
Straightforward calculations show that $\mu_+\equiv\mu_-\equiv\mu_0+v \pmod 2$. Because $H_1(X_{\alpha},V;\R)=0$, the chain complex $\mathcal{H}_{\alpha}$ has two portions: 
$0\rightarrow H_0(V;\R)\rightarrow H_0(
X_{\alpha};\R)\rightarrow H_0(X_{\alpha},V;\R)\rightarrow 0$, and
\begin{multline*}0\rightarrow H_2(V;\R)\rightarrow H_2(X_{\alpha};\R)\rightarrow H_2(X_{\alpha},V;\R)\rightarrow H_1(V;\R)
\rightarrow H_1(X_{\alpha};\R)\rightarrow 0.
\end{multline*}
For the purpose of comparison we only need to look at the second portion.

 When $\alpha=+$: Recalling that $\dim(H_2(V;\R))$ is the same as $\dim(H_2(X_{\alpha};\R))$, we see that the torsion of $\mathcal{H}_+$ is the torsion of the chain  $0\rightarrow H_2(X_{\alpha},V;\R)\stackrel{\partial}\rightarrow H_1(V;\R)
\rightarrow H_1(X_{\alpha};\R)\rightarrow 0$. Since  $[\partial D_+]=[bd^{-1}]=t_v-t_{v+1}$, it follows that $\tau(\mathcal{H}_+)$ is the determinant of the change of bases matrix $[(t_v-t_{v+1},t_1,\dotsc,t_v)/(t_1,\dotsc,t_v,t_{v+1})]$, which is $(-1)^{v+1}$.

When $\alpha=-$: In this case $[\partial D_-]=[a^{-1}c]=t_{v+1}-t_v$, thus $\tau(\mathcal{H}_-)=(-1)^v$.

When $\alpha=0$: The torsion $\tau(\mathcal{H}_0)$ is the torsion of the chain $0\rightarrow H_2(V;\R)\stackrel{i_*}\rightarrow H_2(X_0;\R)\stackrel{j_*}\rightarrow H_2(X_0,V;\R)\rightarrow 0$. The map $i_*$ is an injection,
$i_*(q_i)=q_i,\ 1\le i\le v-1$. The disk $D_0$ is a representative of a generator of $H_2(X_0,V;\R)$. We need to take a lift of $[D_0]$ under the map $j_*$. The union of $D_0$  with part of the boundary of $V$ constitutes either one of the two boundary components of $X_0$ corresponding to $q_v$ and $q_{v+1}$. See Fig. \ref{1component3}, in which the solid two holes torus contains the ball $B$ and the disk $D_0$, while $V$ is outside.
\begin{figure}[h]
\begin{center}
\includegraphics[height=0.1\textheight]{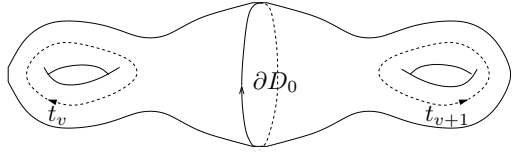}
\end{center}
\caption{The disk $D_0$.}
\label{1component3}
\end{figure}
Because of the chosen orientation of $\partial D_0$ the two corresponding elements in $H_2(X_0)$, which are lifts of $[D_0]$ under $j_*$, are $-q_v$ and $q_{v+1}=-(q_1+q_2+\dotsb+q_v)$. The choice of either lift would result that $\tau(\mathcal{H}_0)=[(q_1,\dotsc,q_{v-1},-q_v)/(q_1,\dotsc,q_v)]=-1$.

Collecting the above computations and comparisons of $\mu_{\alpha}$ and $\tau(\mathcal{H}_{\alpha})$ we conclude that $\breve{\tau}(C_*(X_+,\R))=-\breve{\tau}(C_*(X_-,\R))=\breve{\tau}(C_*(X_0,\R))$. 

{\it Case 2: The two strands of $L_+$ at the crossing come from  different
  components.} See Fig. \ref{2components}. 
\begin{figure}[h]
\begin{center}
\includegraphics[height=0.17\textheight]{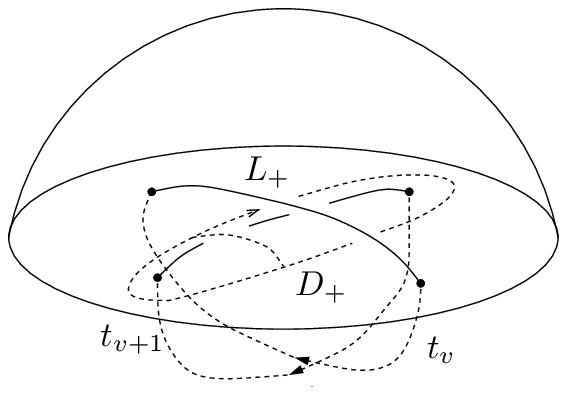}
\end{center}
\caption{Case 2.}
\label{2components}
\end{figure}
Similar to Case 1, the comparison of  $\breve\tau(C_*(X_{\alpha};\R))$ is reduced to the comparison of $(-1)^{\mu_{\alpha}}\text{sign}(\tau(\mathcal{H}_{\alpha}))$. Straightforward calculations  give that $\mu_+\equiv\mu_-\equiv\mu_0+v\pmod 2$. Again to study  $\tau(\mathcal{H}_{\alpha})$ we only need to pay attention to the exact chain complex
\begin{multline*}0\rightarrow H_2(V;\R)\rightarrow H_2(X_{\alpha};\R)\rightarrow H_2(X_{\alpha},V;\R)\rightarrow H_1(V;\R)
\rightarrow H_1(X_{\alpha};\R)\rightarrow 0.
\end{multline*}

When $\alpha=+$: $\tau(\mathcal{H}_+)$ is the torsion of the chain $0\rightarrow H_2(V;\R)\rightarrow H_2(X_+;\R)\rightarrow H_2(X_+,V;\R)\rightarrow 0$. The lift of $[D_+]\in H_2(X_+,V;\R)$ to $H_2(X_+;\R)$ is either $q_v$ or $-q_{v+1}$. With either lift the we have $\tau(\mathcal{H}_+)=[(q_1,\dotsc,q_{v-1},q_v)/(q_1,\dotsc,q_v)]=1$.

When $\alpha=-$: Just as the case $\alpha=+$, except that now the lift of $[D_-]$ can be either $-q_v$ or $q_{v+1}$, so $\tau(\mathcal{H}_-)=-1$.

When $\alpha=0$: $\tau(\mathcal{H}_0)$ is the torsion of the chain $0\rightarrow H_2(X_0,V;\R)\stackrel{\partial}\rightarrow H_1(V;\R)
\rightarrow H_1(X_0;\R)\rightarrow 0.$ Since $[\partial D_0]=[a^{-1}b]=t_{v+1}-t_v\in H_1(V;\R)$ we have $\tau(\mathcal{H}_0)=[(t_{v+1}-t_v,t_1,\dots,t_v)/(t_1,\dots,t_{v+1})]=(-1)^v$.

Thus as in Case 1, $\breve{\tau}(C_*(X_+,\R))=-\breve{\tau}(C_*(X_-,\R))=\breve{\tau}(C_*(X_0,\R))$.

Now Formula (\ref{signrefinedformula}) and the skein relation for unrefined torsion (\ref{skeintorsion}) give us a skein relation for sign-refined torsion:
\begin{equation}\label{skeinrefinedtorsion}\tau_{0,\ L_+}^{\varphi'}(t)+(1-t)\tau_{0,\ L_0}^{\varphi'}(t)-t\tau_{0,\ L_-}^{\varphi'}(t)=0,
\end{equation} 
provided that $L_+$, $L_-$ and $L_0$ belong to the same torsion class. 
\end{subsubsection}
\begin{subsubsection}{Definition of the normalized one variable twisted Alexander function}
For a given link $L$ the sign-refined torsion $\tau_{0,\ L}^{\varphi'}(t)$ is defined up to $t^n,\ n\in\Z$. Using Theorem \ref{symmetrytorsion}, there is a number $r\in\Z$ arising from the symmetry of (un-refined) torsion such that 
$\tau^{\varphi'}_{0,\ L}(t^{-1})=\pm t^r\tau^{\varphi'}_{0,\ L}(t)$  as
elements in $\Q(t)$. 

Define the normalized twisted Alexander function of a link $L$ to be
\begin{equation}\nabla_L(t)=-t^r\tau_{0,\ L}^{\varphi'}(t^2)\label{defnormalized}.
\end{equation} 
\index{twisted Alexander function!normalized}

Notice that 
$\nabla_{L}(t^{-1})=-t^{-r}\tau_{0,\ L}^{\varphi'}(t^{-2})=\pm t^{-r}t^{2r}\tau_{0,\ L}^{\varphi'}(t^2)=\pm
t^r\tau_{0,\ L}^{\varphi'}(t^2)=\pm\nabla_L(t)$. Thus $\nabla_L(t)$ is
symmetric, up to a sign. From Theorems \ref{otAPtorsion} and
\ref{otAPnontorsion}, the function $\nabla_L(t)$ is an element of $\Z[t^{\pm 1}]$ (a Laurent polynomial) if $L$ is nontorsion, and is an element of $\Z[t^{\pm
    1},(t-t^{-1})^{-1}]$ (a Laurent polynomial divided by $(t-t^{-1})^n$) if $L$ is torsion.
    
\begin{proposition}The function $\nabla_L(t)$ does not depend on the
  choice of a representative of $\tau_{0,\ L}^{\varphi'}(t)$ and so is completely defined without indeterminacies.
\end{proposition}

\begin{proof}Suppose that $\tau$ and $\tau'$ are two  representatives
  of the (sign-refined) torsion $\tau^{\varphi'}_{0,\ L}$. Then $\tau'(t)=t^m\tau(t)$ for
  some $m\in\Z$. This
  implies that there is an $n\in\Z$ such that $\nabla'(t)=
  t^n\nabla(t)$. Since $\nabla(t^{-1})=\pm\nabla(t)$ and
  $\nabla'(t^{-1})=\pm\nabla'(t)$ we must have $n=0$, that is
  $\nabla'(t)=\nabla(t)$.
\end{proof}

\end{subsubsection}

\end{subsection}

\subsection{A skein relation for the normalized twisted Alexander function}
\begin{theorem}\label{mainthmc1}If $L_+$, $L_-$ and $L_0$ belong to the same torsion class then the normalized one variable twisted Alexander
  function satisfies the skein relation:
\begin{equation}\label{skeinfinal}\nabla_{L+}(t)-\nabla_{L_-}(t)=(t-t^{-1})\nabla_{L_0}(t).
\end{equation}
\end{theorem}

\begin{proof}
Replacing $t$ by $t^2$ in Eq.  (\ref{skeinrefinedtorsion}), and using
Eq. (\ref{defnormalized}) we have
$$t^{-r_+}\nabla_{L_+}(t)+(1-t^2)t^{-r_0}\nabla_{L_0}(t)-t^{2-r_-}\nabla_{L_-}(t)=0,$$
that is
$$\nabla_{L_+}(t)=(t-t^{-1})t^{1+r_+-r_0}\nabla_{L_0}(t)+
  t^{2+r_+-r_-}\nabla_{L_-}(t).$$

Let $u=2+r_+-r_-$ and $v=1+r_+-r_0$ we get
\begin{equation}\label{skeinnabla}\nabla_{L_+}(t)=(t-t^{-1})t^v\nabla_{L_0}(t)+
  t^u\nabla_{L_-}(t).
\end{equation}

The purpose of the rest of the proof is to show that $u=v=0$. The
idea is to show that $u$ and $v$ are independent
of the link. This is achieved by studying the numbers
$r_{\alpha}$'s. Since these numbers arise from the symmetry of torsion, a
study of duality of torsion is needed. 

Topologically the complement $X_{\alpha}$ of $L_{\alpha}$ is the union
of $V$ and a $2$-handle $H_{\alpha}$ glued to $V$ along the loop $\partial
D_{\alpha}$. Assume that $X_{\alpha}$ is triangulated by a triangulation of $V$ together with a
compatible triangulation of $H_{\alpha}$. Let $\widetilde{X}_{\alpha}$ be the $D=\Z\times\Tors H_1(X_{\alpha})$ cover of $X_{\alpha}$ corresponding to the kernel of the map $proj_{\alpha}:\pi_1(X_{\alpha})\rightarrow H_1(X_{\alpha})\rightarrow G\times\Tors H_1(X_{\alpha})\rightarrow \{t^m:m\in\Z\}\times\Tors H_1(X_{\alpha})$. As in Section \ref{C-Calpha}, $\widetilde{X}_{\alpha}$ can be constructed as 
$\widetilde{V}\cup_{t\in D} t\widetilde{H}_{\alpha}$, i.e. $\widetilde{V}$ with disjoint copies of $H_{\alpha}$ glued in along the lifts of $\partial D_{\alpha}$.
Because of our assumption that  $L_{\alpha}$'s belong to the same torsion class, the deck transformation group $D$ does not depend on $\alpha$. An induced triangulation $Y_{\alpha}$
of $\widetilde{X}_{\alpha}$ is obtained, which is equivariant under the
action of $D$. Let
$Y^*_{\alpha}$ be its dual cell decomposition and $\partial
Y^*_{\alpha}$ be the restriction of $Y^*_{\alpha}$ to the boundary
$\partial \widetilde{X}_{\alpha}$.

Let $E_{\alpha}=\Q(t)\otimes_{\varphi'} C_*(Y_{\alpha})$,
$F_{\alpha}=\Q(t)\otimes_{\varphi'} C_*(Y^*_{\alpha})$, 
$\partial F_{\alpha}=\Q(t)\otimes_{\varphi'} C_*(\partial
Y^*_{\alpha})$.

Choose a fundamental family of cells $e_{\alpha}$ for
$Y_{\alpha}$ such that all the cells in $e_{\alpha}$ that cover a cell in $H_{\alpha}$ are contained in the same $\widetilde{H}_{\alpha}$. Denote by $e^*_{\alpha}$ the family of cells in $Y^*_{\alpha}$ that are dual to the simplexes in $e_{\alpha}$.

The proof consists of the following steps.

{\it Step 1: Studying $\tau(F_{\alpha})$.} The triangulation $Y_{\alpha}$ and its  dual cell decomposition $Y^*_{\alpha}$ has a common
cellular subdivision, namely the first barycentric subdivision $Y'_{\alpha}$ of
$Y_{\alpha}$. It is possible to choose two fundamental family of cells
for $\widetilde{X}_{\alpha}$ corresponding to $Y'_{\alpha}$. The first is $a$, consisting of the cells
$a_1,a_2,\dotsc,a_n$, each of which is contained in a cell in
$e_{\alpha}$. This provides a chosen basis for $E_{\alpha}$. The second fundamental family of cells is
$b$, consisting of the cells $b_1,b_2,\dotsc,b_n$, each of which is
contained in a cell in $e^*_{\alpha}$, providing a chosen basis for  $F_{\alpha}$. 

Using invariance of torsion under cellular subdivision (see \cite[Lemma 4.3.3 iii]{Turaev86}) we have
$\tau(E_{\alpha},e_{\alpha})=\pm\tau(\Q(t)\otimes_{\varphi'}
C_*(Y'_{\alpha}),a)$ and $\tau(F_{\alpha},e^*_{\alpha})=\pm\tau(\Q(t)\otimes_{\varphi'}
C_*(Y'_{\alpha}),b)$. Let us compare the torsion of the same chain
complex $\Q(t)\otimes_{\varphi'}
C_*(Y'_{\alpha})$ with different bases $a$ and $b$. 

We have $\tau(\Q(t)\otimes_{\varphi'}
C_*(Y'_{\alpha}),b)=\tau(\Q(t)\otimes_{\varphi'}
C_*(Y'_{\alpha}),a)\varphi'([b/a])$, where $[b/a]\in D$ denotes
the determinant of the change of base matrix.  
If two cells $a_i$ and  $b_j$ cover the same cell in the $2$-handle
$H_{\alpha}$ then they must be contained in the same $\widetilde{H}_{\alpha}$ because of our choice for $e_{\alpha}$ above, and so 
$a_i$ and $b_j$ must be the
same cell. This means that the correctional term $\varphi'([b/a])$ does not
depend on $\alpha$. 

Thus there is 
$\beta\in\Z$ which does not depend on $\alpha$ such that 
\begin{equation}\label{torsionFbeta}
\tau(F_{\alpha},e^*_{\alpha})=\pm
t^{\beta}\tau(E_{\alpha},e_{\alpha})=\pm
t^{\beta}\tau_{L_{\alpha}}^{\varphi'}(t).
\end{equation}

{\it Step 2: Studying the chain $\partial F_{\alpha}$.} Consider the short exact sequence of chain complexes
\begin{equation}\label{seqF}
0\rightarrow \partial F_{\alpha}\rightarrow F_{\alpha}\rightarrow
F_{\alpha}/\partial F_{\alpha}\rightarrow 0.
\end{equation}
Note that $\partial X_{\alpha}$ is a collection of tori. It is simple to see that the chain $\partial F_{\alpha}$ is exact and its
  torsion -- the  torsion of a collection of tori -- is $1$ up to $\pm
  t^n$.

The long homological exact sequence associated with the short exact
sequence (\ref{seqF}) above shows that $F_{\alpha}$ is exact if and only
if $F_{\alpha}/\partial F_{\alpha}$ is exact. Note that by the
invariance of torsion under cellular subdivisions, $F_{\alpha}$
is exact if and only if $E_{\alpha}$ is exact, and in any case
$\tau(F_{\alpha})=\tau(E_{\alpha})$ up to $\pm t^n$. The
product formula for torsion of chain complexes applied to the short exact sequence
(\ref{seqF}) gives 
\begin{equation}\label{productF}
\tau(F_{\alpha})=\pm\tau(\partial
F_{\alpha})\tau(F_{\alpha}/\partial F_{\alpha}). 
\end{equation}
Both sides are zero
when $F_{\alpha}$ is not exact.

Let $R$ be the union
of those tori of $\partial X_{\alpha}$
which do not involve the crossing, i.e. $R\cap B=\emptyset$, where
$B$ is the ball enclosing the crossing under scrutiny as in Fig. \ref{skein2}. Then
$\partial X_{\alpha}\setminus R$ is a disjoint union of two tori if
the two strands at the crossing belong to different components of the
link or it is
just a torus if the two strands belong to the same component.

Let $P=\Q(t)\otimes_{\varphi'}C_*(\partial Y^*|_R)$ and
$Q_{\alpha}=\Q(t)\otimes_{\varphi'}C_*(\partial Y^*|_{\partial
  X_{\alpha}\setminus R})$. Then $\partial F_{\alpha}=P\oplus
Q_{\alpha}$. Note that $\partial F_{\alpha}$, $P$ and $Q_{\alpha}$
are all acyclic chain complexes. The torsion of $P$ does not depend on
$\alpha$ and is $1$ up to units: $\tau(P)= \pm t^p$ for some
$p\in\Z$, on the other hand $\tau(Q_{\alpha})=\pm t^{q_{\alpha}}$ for
some $q_{\alpha}\in\Z$. The number $q_{\alpha}$ depends on how the
lifting cells are chosen. It depends  only on whether the two strands at the crossing  under
investigation belong to the same component or two different components
of the link  $L_{\alpha}$. The product formula gives us
$\tau(\partial F_{\alpha})=\pm\tau(P)\tau(Q_{\alpha})=\pm t^{p+{q_{\alpha}}}$.  
 
{\it Step 3: Studying $\tau(F_{\alpha}/\partial F_{\alpha})$.} 
By the symmetry of torsion (Section \ref{symmetry}), $\tau(F_{\alpha}/\partial
F_{\alpha})=\overline{\tau(E_{\alpha})}=\tau_{L_{\alpha}}^{\varphi'}(t^{-1}) = \pm
t^{r_{\alpha}}\tau_{L_{\alpha}}^{\varphi'}(t)$. Note that this
$r_{\alpha}$ is the one in Eq. (\ref{defnormalized}).

{\it Step 4:  Skein relation for $\nabla$.} From Eq. (\ref{productF}), Step 2 and Step 3 we have $\tau(F_{\alpha})=\pm
t^{p+{q_{\alpha}}}
t^{r_{\alpha}}\tau_{L_{\alpha}}^{\varphi'}(t)$. Comparing with Eq.
(\ref{torsionFbeta}) we get $\pm
t^{p+{q_{\alpha}}+r_{\alpha}}\tau_{L_{\alpha}}^{\varphi'}(t)=\pm
t^{\beta}\tau_{L_{\alpha}}^{\varphi'}(t)$. This 
gives us 
\begin{equation}\label{betapqr}
\beta=\beta_{\alpha}=p+{q_{\alpha}}+r_{\alpha}.
\end{equation}
Using Eq. (\ref{betapqr})  we have $u=2+r_+-r_-=2+q_--q_+$  and
$v=1+r_+-r_0=1+q_0-q_+$. Thus Eq. (\ref{skeinnabla}) depends on the links
$L_{\alpha}$'s only to the extent that whether the two strands at the crossing  under
investigation belong to the same component or two different components
of the link  $L_{\alpha}$. Equation (\ref{skeinnabla})
is satisfied with the same $u$ and $v$ for all link $L_+$ whose two
strands at the crossing come from the same component, and is also
satisfied with the same $u$ and $v$ for all link $L_+$ whose two
strands at the crossing come from two different components. Thus in
each case a particular example is enough to determine the values of
$u$ and $v$.

{\it Case 1: The two strands of
$L_+$ at the
crossing come from one component.} Consider the knot $3_1$ and the particular crossing in
Fig. \ref{3_1-5_6}. 
\begin{figure}[h]
\begin{center}
\includegraphics[height=0.13\textheight]{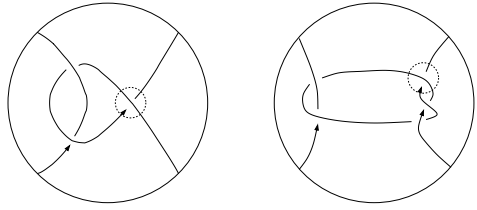}
\end{center}
\caption{The knots $3_1$ and $5_6$.}
\label{3_1-5_6}
\end{figure}
Direct computation gives that $\nabla_{L_+}(t)=\pm(t-t^{-1})$,
$\nabla_{L_-}(t)=\pm(t-t^{-1})$, and $\nabla_{L_0}(t)=0$, thus $u=0$. Also consider the knot $5_6$ in that figure. We have $\nabla_{L_+}(t)=\pm(t-t^{-1})$, $\nabla_{L_-}(t)=\pm(t-t^{-1})(t^2-1+t^{-2})$, and $\nabla_{L_0}(t)=\pm (t-t^{-1})^2$, thus $v=0$.

{\it Case 2: The two strands of
$L_+$ at the
crossing come from  two different components.} Consider the link
$4^2_2$ in Fig. \ref{4^2_2-link}. 
\begin{figure}[h]
\begin{center}
\includegraphics[height=0.13\textheight]{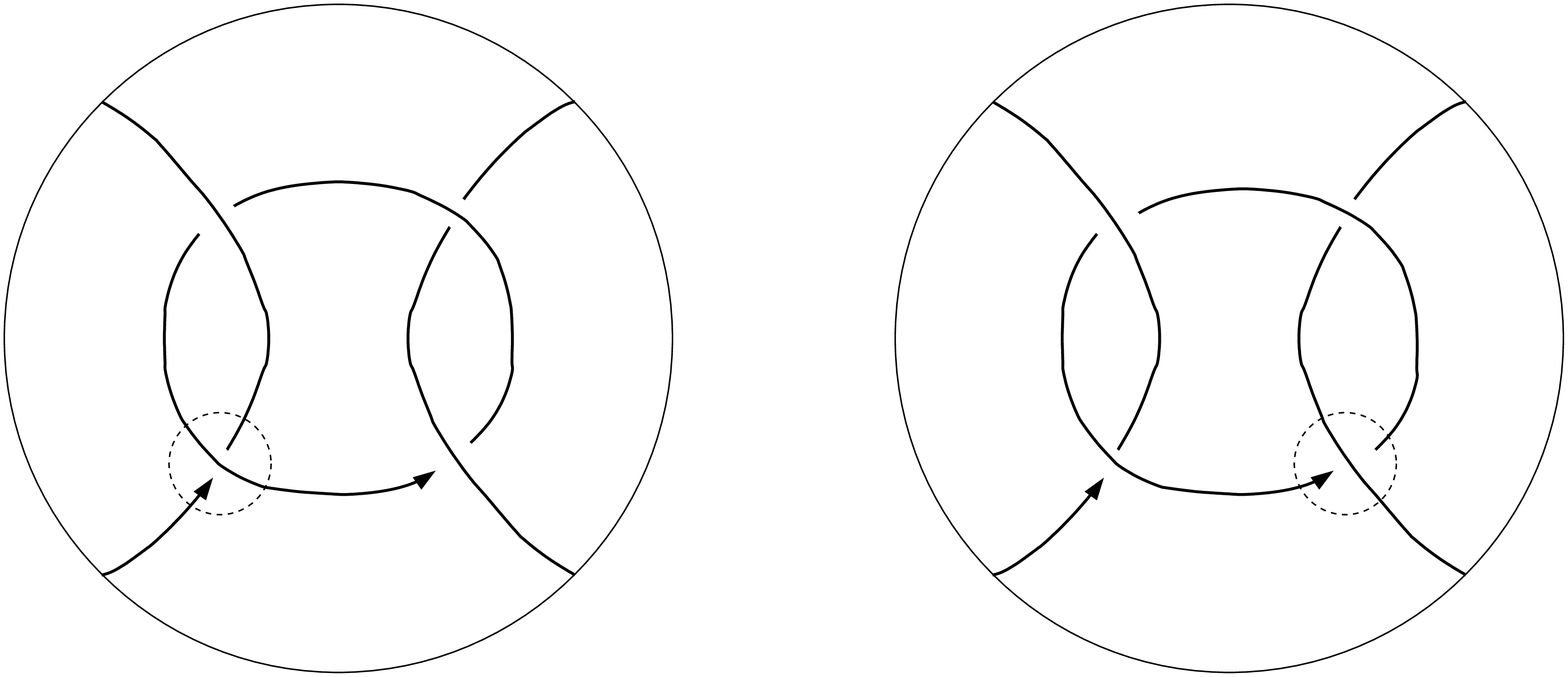}
\end{center}
\caption{The link $4^2_2$.}
\label{4^2_2-link}
\end{figure}
At the first crossing in the
figure, $\nabla_{L_+}(t)=\pm(t-t^{-1})^2$,
$\nabla_{L_0}(t)=\pm(t-t^{-1})$, and $\nabla_{L_-}(t)=0$, thus $v=0$. On the other hand at the second crossing in the figure $\nabla_{L_-}(t)=\pm(t-t^{-1})^2$,
$\nabla_{L_0}(t)=\pm(t-t^{-1})$, and $\nabla_{L_+}(t)=0$, thus $u=0$.

In both cases $u=v=0$, and the proof of Theorem \ref{mainthmc1} is completed. 
\end{proof}

\begin{remark}In general it is not possible to compute $\nabla_{L}(t)$ from the skein relation (\ref{skeinfinal}) alone because of the restriction of our theorem that the torsion classes do not change after a smoothing at a crossing.
\end{remark}

\end{section}

\section{Relationships among twisted and untwisted Alexander polynomials}
Suppose that $L$ is a nontorsion link \index{link!nontorsion} in $\RP^3$. Let $\widetilde{L}$ be the preimage of $L$ under the canonical covering map from $S^3$ to $\RP^3$. Because each component of $L$ is null-homologous hence is null-homotopic in $\RP^3$, its preimage in $S^3$ has two components. Thus $\widetilde{L}$ has an even number of components. 
A way to draw a diagram for $\widetilde{L}$ is to put a copy of a diagram $D$ of $L$ on the top disk of a cylinder. On the bottom disk put a diagram obtained from $D$ by reflecting it through the center of the disk, then connect the corresponding boundary points on the boundary circles of the top and bottom disks by vertical lines. If furthermore we rotate the bottom disk an angle of $180^{\circ}$ along a horizontal line passing through the disk center (i.e. flipping it, changing undercrossings to overcrossings and vice versa) then we obtain Drobotukhina's description in \cite[p. 616]{Drobotukhina90}.


Consider the following diagram of coverings:
\begin{equation}\label{diagramofcoverings}
\xymatrix{& {\widetilde{X}}\ar[ld]_{p_2}^{\Z_2}\ar[rd]^{p_4}_G\ar[dd]_p^{G\times\Z_2}&\\
{\widetilde{X}_G} \ar[rd]^{p_1}_{G}&  &{\widetilde{X}_2=S^3\setminus \stackrel{\circ}{N(\widetilde{L})}}\ar[ld]_{p_3}^{\Z_2}\\
& X=\RP^3\setminus \stackrel{\circ}{N(L)}&}
\end{equation}
In the diagram $p:\widetilde{X}\rightarrow X$ corresponds to the kernel of the map $\pi\rightarrow H$; $p_1:\widetilde{X}_G\rightarrow X$ corresponds to the kernel of the map $\pi\rightarrow H\rightarrow G$; $p_3:\widetilde{X}_2\rightarrow X$ corresponds to the kernel of the map $\pi\rightarrow H\rightarrow \Z_2$; and $p_2$ and $p_4$ are lifts of $p$. The diagram is commutative. The cellular structure of $X$ induces cellular structures on the remaining spaces.

 Let $C_i^+(\widetilde{X})$ be the subcomplex of $C_i(\widetilde{X})$ generated by chains of the form $\sigma+u\sigma$ where $\sigma$ is an $i$-cell in $\widetilde{X}$. Similarly let $C_i^-(\widetilde{X})$ be the subcomplex generated by chains of the form $\sigma-u\sigma$. Consider $\Q(G)\otimes_{\Z[H],\varphi}C_i(\widetilde{X})$, where $\varphi$ is the twisted map of Section \ref{twistedmap}.

\begin{proposition}We have the following isomorphisms of $\Q(G)$-vector spaces: 

a). $\Q(G)\otimes_{\Z[H]}C_i(\widetilde{X})=(\Q(G)\otimes_{\Z[H]}C_i^+(\widetilde{X}))\oplus (\Q(G)\otimes_{\Z[H]}C_i^-(\widetilde{X}))$.

b). $\Q(G)\otimes_{\Z[H]}C_i(\widetilde{X}_G)\cong \Q(G)\otimes_{\Z[H]}C_i^+(\widetilde{X})$.

c). $\Q(G)\otimes_{\Z[H],\varphi}C_i(\widetilde{X})\cong \Q(G)\otimes_{\Z[H]}C_i^-(\widetilde{X})$.

\end{proposition}
\begin{proof}Here we are dealing with homology with local coefficients\index{local coefficients!homology} and the following proof is adapted from Hatcher \cite[p. 330]{Hatcher01}. 

a).  Noting that $C_i^+(\widetilde{X})\cap C_i^-(\widetilde{X})=\{0\}$ and $\sigma=((\sigma+u\sigma)+(\sigma-u\sigma))/2$, the result follows immediately.

b). A cell in $\widetilde{X}$ is a lift of a cell in $\widetilde{X}_G$. 
The isomorphism is induced from the map $\sigma\mapsto (\tilde{\sigma}+u\tilde{\sigma})$.

c). Consider the the projection $pr:\Q(G)\otimes_{\Z[H]}C_i(\widetilde{X})\rightarrow \Q(G)\otimes_{\Z[H],\varphi}C_i(\widetilde{X})$ mapping $1\otimes\sigma$ to $1\otimes_{\varphi}\sigma$. We have $pr(1\otimes(\sigma+u\sigma))=1\otimes_{\varphi}\sigma+1\otimes_{\varphi}u\sigma=1\otimes_{\varphi}\sigma+1\cdot u\otimes_{\varphi}\sigma=1\otimes_{\varphi}\sigma+\varphi(u)\otimes_{\varphi}\sigma=0$, since $\varphi(u)=-1$. This implies $\Q(G)\otimes_{\Z[H]}C_i^+(\widetilde{X})\subset \ker(pr)$. By a similar argument we see that $(\Q(G)\otimes_{\Z[H]}C_i^-(\widetilde{X}))\cap \ker(pr)=\{0\}$. Thus $\Q(G)\otimes_{\Z[H]}C_i^+(\widetilde{X})= \ker(pr)$ and using a) the result follows.
\end{proof}

Note that $\Q(G)\otimes_{\Z[H]}C_i(\widetilde{X}_G)$ is in fact $\Q(G)\otimes_{\Z[G]}C_i(\widetilde{X}_G)$. It follows from this proposition that we have the short exact sequence of chain complexes of $\Q(G)$-vector spaces:
\begin{equation}\label{seqAtA}0\rightarrow \Q(G)\otimes_{\Z[H],\varphi}C_*(\widetilde{X})\rightarrow \Q(G)\otimes_{\Z[H]}C_*(\widetilde{X})\rightarrow \Q(G)\otimes_{\Z[G]}C_*(\widetilde{X}_G)\rightarrow 0.
\end{equation}

From this sequence we now derive a relationship among multi-variable Alexander polynomials. If $L$ is a nontorsion link having $v$ components then $\widetilde{L}$ has $2v$ components. We enumerate so that the $i$-th component and the $(v+i)$-th component of $\widetilde{L}$ are projected  to the same $i$-th component of $L$. Let $\psi$ be the homomorphism from $\Z[t_1^{\pm 1} ,t_2^{\pm 1},\dotsc,t_v^{\pm 1},t_{v+1}^{\pm 1},\dotsc,t_{2v}^{\pm 1}]$ to 
$\Z[t_1^{\pm 1} ,t_2^{\pm 1},\dotsc,t_v^{\pm 1}]$ identifying $t_{v+i}$ with $t_i$ for all $1\le i\le v$. Consider the multi-variable Alexander polynomial of $\widetilde{L}$, $\Delta_{\widetilde{L}}(t_1,t_2,\dotsc,t_{2v})$. Let $\Delta'_{\widetilde{L}}(t_1,t_2,\dotsc,t_v)$ be obtained from $\Delta_{\widetilde{L}}(t_1,t_2,\dotsc,t_{2v})$ by identifying the $t_i$ and $t_{v+i}$ variables for all $1\le i\le v$, that is $\Delta'(\widetilde{L})=\psi(\Delta(\widetilde{L}))$. Recall from our fixed splitting of $H$ in Section \ref{twistedmap} that the free part $G$ is generated by the meridians of the components of $L$, thus $\Z[G]=\Z[t_1^{\pm 1} ,t_2^{\pm 1},\dotsc,t_v^{\pm 1}]$.

\begin{example}Let $K$ be the knot $2_1$ in Drobotukhina's table (see Example \ref{2_1-knotex}). Then $\Delta_K^{\varphi}(t)=t-1$ and $\Delta_K(t)=t^2+1$. The lift $\widetilde{K}$ of this knot is the link $4^1_2$ in Rolfsen's table,  $\Delta_{\widetilde{K}}(t_1,t_2)=t_1t_2+1$, and $\Delta'_{\widetilde{K}}(t)=t^2+1$.
\end{example}

\begin{theorem}\label{cover} Let $L$ be a nontorsion link. If $L$ has one component then $(t-1)\Delta'(\widetilde{L})=\Delta(L)\Delta^{\varphi}(L)$ as elements in $\Z[t^{\pm 1}]$. If $L$ has at least two components then $\Delta'(\widetilde{L})=\Delta(L)\Delta^{\varphi}(L)$ as elements in $\Z[t_1^{\pm 1} ,t_2^{\pm 1},\dotsc,t_v^{\pm 1}]$.
\end{theorem}

\begin{proof}Recall the diagram of covering spaces (\ref{diagramofcoverings}). The map $p_4$ corresponds to the kernel of the canonical projection $\pi_1(\widetilde{X}_2)\rightarrow H_1(\widetilde{X}_2)=\langle t_1,\dotsc,t_{2v}/t_it_j=t_jt_i\rangle\rightarrow G=\langle t_1,\dotsc,t_v/t_it_j=t_jt_i\rangle$, where the second projection identifies $t_i$ and $t_{v+i}$ for all $1\le i\le v$. Thus $p_{4*}(\pi_1(\widetilde{X}))$ will be the subgroup of $\pi_1(\widetilde{X}_2)$ whose projection to $H_1(\widetilde{X}_2)$ is $\{t_1^{\alpha_1}\dotsb t_{2v}^{\alpha_{2v}}/\alpha_i+\alpha_{v+i}=0,1\le i\le v\}$. Then $p_{3*}$ will send $p_{4*}(\pi_1(\widetilde{X}))$ to the subgroup of $\pi_1(X)$ whose projection to $H$ is $\{t_1^{\alpha_1+\alpha_{v+1}}\dotsb t_v^{\alpha_v+\alpha_{2v}}\}=\{1\}$. So $(p_3\circ p_4)_*$ sends $\pi_1(\widetilde{X})$ to the subgroup of $\pi$ which vanishes in $H$, this is why $p_3\circ p_4=p$.

Now we look at the space $\widetilde{X}$ as the $G$-cover of $\widetilde{X}_2$ corresponding to $p_4$. Then there is an action of $G$ on $C_i(\widetilde{X})$ turning it to a $\Z[G]$-module $C_i'(\widetilde{X})$, and so we can form the vector space $\Q(G)\otimes_{\Z[G]}C_i'(\widetilde{X})$. It can be seen that $\Q(G)\otimes_{\Z[G]}C_i'(\widetilde{X})\cong \Q(G)\otimes_{\Z[H]}C_*(\widetilde{X})$. Thus the sequence (\ref{seqAtA}) becomes
$$0\rightarrow \Q(G)\otimes_{\Z[H],\varphi}C_*(\widetilde{X})\rightarrow\Q(G)\otimes_{\Z[G]}C_*'(\widetilde{X})\rightarrow \Q(G)\otimes_{\Z[G]}C_*(\widetilde{X}_G)\rightarrow 0.$$
Apply the product formula for torsion (\ref{productformulasforun-refinedtorsion}) to this short exact sequence we obtain 
\begin{equation}\label{finalexactseq}\tau(\Q(G)\otimes_{\Z[G]}C_*'(\widetilde{X}))=\pm\tau(\Q(G)\otimes_{\Z[G]}C_*(\widetilde{X}_G))\tau(\Q(G)\otimes_{\Z[H],\varphi}C_*(\widetilde{X})).
\end{equation}

Theorem \ref{torsion=tAP} says that $\tau(\Q(G)\otimes_{\Z[H],\varphi}C_*(\widetilde{X}))$ is $\Delta^{\varphi}(L)$; Remark \ref{Ap=R} says $\tau(\Q(G)\otimes_{\Z[G]}C_*(\widetilde{X}_G))$ is $\Delta(L)$ if $L$ has more than one component and is $\Delta(L)/(t-1)$ if $L$ has one component. Finally the identification of torsion and Alexander polynomial for links in $S^3$ (\cite{Milnor62}, \cite[p. 55]{Turaev01}) says that $\tau(\Q(G)\otimes_{\Z[G]}C_*'(\widetilde{X}))$ is $\psi(\Delta(\widetilde{L}))$ if $\widetilde{L}$ has more than one component (here the functority of torsion \cite[Lemma 13.5]{Turaev01} is used). The theorem then follows from (\ref{finalexactseq}).
\end{proof}

\begin{remark}A similar result also holds true if we consider only one variable polynomials.
\end{remark}

\section*{Acknowledgments}We wish to thank Vladimir Turaev for raising many of the questions considered in this paper. We thank the referee for valuable comments. One of the authors (V.H.) would like to thank the Department of Mathematics at the State University of New York at Buffalo, where part of this work was carried out.

\bibliographystyle{amsalpha}
\bibliography{my}

\end{document}